\documentclass[12pt]{amsart}

\usepackage{CJKutf8}
\usepackage{amssymb}
\usepackage{eucal}
\usepackage{color}

\usepackage{multirow,bigdelim}
\usepackage{amsfonts}
\usepackage{dsfont}
\usepackage{amsmath}
\usepackage{mathrsfs}
\usepackage{todonotes}
\usepackage{hyperref}
\usepackage{multimedia}
\usepackage{graphicx}
\usepackage{indentfirst}
\usepackage{bm}
\usepackage{amsthm}
\usepackage{mathrsfs}
\usepackage{latexsym}
\usepackage{amsmath}
\usepackage{amssymb}
\usepackage{microtype}
\usepackage{relsize}
\usepackage{hyperref}

\allowdisplaybreaks[4]
\DeclareSymbolFont{SY}{U}{psy}{m}{n}
\DeclareMathSymbol{\emptyset}{\mathord}{SY}{'306}

\theoremstyle{plain}

\newtheorem{thm}{Theorem}[section]
\newtheorem{cor}[thm]{Corollary}
\newtheorem{lem}[thm]{Lemma}
\newtheorem{prop}[thm]{Proposition}
\theoremstyle{definition}
\newtheorem{defn}[thm]{Definition}
\newtheorem{rem}[thm]{Remark}
\newtheorem{ex}[thm]{Example}

\numberwithin{equation}{section}

\def\H{{\mathcal H}}

\def\beq{\begin{eqnarray}}
\def\eeq{\end{eqnarray}}
\def\beqa{\begin{eqnarray*}}
\def\eeqa{\end{eqnarray*}}

\begin{document}
\title[Geometry of holomorphic vector bundles and similarity]{Geometry of holomorphic vector bundles and similarity of commuting tuples of operators}

\author{Yingli Hou, Kui Ji, Shanshan Ji$^*$ and Jing Xu}
\curraddr[Y. Hou, K. Ji, S. Ji and J. Xu ]{School of Mathematical Sciences, Hebei Normal University,
Shijiazhuang, Hebei 050016, China}

\email[Y. Hou]{houyingli0912@sina.com}
\email[K. Ji]{jikui@hebtu.edu.cn, jikuikui@163.com}
\email[S. Ji]{jishanshan15@outlook.com}
\email[J. Xu]{xujingmath@outlook.com}

\thanks{The first author was supported by National
Natural Science Foundation of China (Grant No. 12001159).
The second author was supported by National Natural Science Foundation of China (Grant No. 11831006 and
11922108).}
\thanks{* Corresponding author}

\begin{abstract}
In this paper, a new criterion for the similarity of commuting tuples of operators on Hilbert spaces is introduced.
As an application, we obtain a geometric similarity invariant of tuples in the Cowen-Douglas class which gives a partial answer to a question
raised by R.G. Douglas in \cite{CD,Douglas} about the similarity of quasi-free Hilbert modules.
Moreover, a new subclass of commuting tuples of Cowen-Douglas class is obtained.


\end{abstract}

\maketitle

\section{Introduction}

Let $\mathcal H$ be a complex separable Hilbert space, and ${\mathcal L}({\mathcal H})$ the collection of bounded linear operators on $\mathcal H$.
Problems of operator theory often involve unitary and similarity equivalences of operators (operator tuples).
For a positive integer $m$, let $\textbf{T}=(T_{1}, \cdots, T_{m})$ be an $m$-tuple of bounded operators acting on $\mathcal{H}$. If $\textbf{T}$ satisfies $T_{i}T_{j}=T_{j}T_{i}$ for all $1\leq i,j\leq m$, $\textbf{T}$ will be referred to as a commuting tuple.
Let $\textbf{S}=(S_{1}, \cdots, S_{m})$ be a commuting tuple and $S_{i}\in \mathcal{L}(\widetilde{\mathcal{H}}), 1\leq i\leq m$.
The equation $X\textbf{T}=\textbf{S}X$ for some $X\in \mathcal{L}(\mathcal{H}, \widetilde{\mathcal{H}})$ means $XT_{i}=S_{i}X$ for all $1\leq i\leq m$.
If $X$ is a unitary operator, $\textbf{T}$ is unitarily equivalent to $\textbf{S}$ (denoted by $\textbf{T}\sim_{u}\textbf{S}$).
If $X$ is invertible, $\textbf{T}$ is similar to $\textbf{S}$ ($\textbf{T}\sim_{s}\textbf{S}$).




Let $\Omega$ be a connected open subset of $\mathbb{C}$ and $n$ be a positive integer.
In \cite{CD}, M.J. Cowen and R.G. Douglas introduced a class of bounded linear operators, denoted by $\mathcal{B}^{1}_n(\Omega)$,
which contains $\Omega$ as eigenvalues of constant multiplicity $n$.
In \cite{CD3}, they also pointed out that some results of \cite{CD} can be directly extended to $\Omega\subset \mathbb{C}^m,m>1$,
i.e. the class of operators $\mathcal{B}^{1}_n(\Omega)$ may be generalized to an operator tuple class $\mathcal{B}^{m}_n(\Omega)$.
Each Cowen-Douglas tuple naturally determines a Hermitian holomorphic vector bundle,
and such two tuples are unitarily equivalent if and only if there is an isometric and connection-preserving bundle map between the bundles \cite{CD,CD3}.
In particular, the unitary classification of tuples in $\mathcal{B}^{m}_1(\Omega)$ involves only the curvature of Hermitian holomorphic bundles.

In \cite{Misra}, G. Misra has introduced and discussed homogeneous operators in $\mathcal{B}^{1}_1(\mathbb{D})$.
By using the curvature as the invariant, these homogeneous operators have been completely characterized.
For homogeneous operators in $\mathcal{B}^{1}_n(\mathbb{D}),n>1$, A. Koranyi and G. Misra analyzed their structure and proved a classification theorem (see \cite{AM}).
In \cite{Misra2}, G. Misra estimated the curvatures of operators in $\mathcal{B}^{1}_1(\Omega)$, and further obtained a widely used curvature inequality,
stating that the curvature $\mathcal{K}_{S^{*}}$ of the backward shift operator dominates the curvature $\mathcal{K}_{T}$ if $T$ is contractive.
Subsequently, G. Misra and N.S.N. Sastry \cite{MS1,MS2} proved that the inequality holds also for curvatures of tuples in $\mathcal{B}^{m}_1(\Omega)$.
Conversely, the fact that the curvature inequality implies that the operator has a stronger contraction than usual case has been proved by S. Biswas, D. K. Keshari and G. Misra in \cite{BKM}.
Other properties of curvature inequality have been discussed in \cite{BG,DMS,DKKS1,DKKS2,MA,WZ}.

In \cite{CS2}, R.E. Curto and N. Salinas linked the above Cowen-Douglas operator theory to the generalized reproducing kernel theory. 
They also discussed the correspondence between the analytic functional Hilbert space with coordinate multiplication $\mathbf{M}_{z}= (M_{z_{1}},\cdots,M_{z_{m}})$ and the canonical module of Cowen-Douglas tuples, proving the following result
\begin{thm}\cite{CS2}
Under mild conditions, the tuple $\text{\bf{M}}_{z}=(M_{z_{1}},\cdots,M_{z_{m}})$ acting on two analytic functional
Hilbert spaces are unitarily equivalent if and only if their normalized reproducing kernel functions are
intertwined by a constant unitary matrix.
\end{thm}

In \cite{C12}, it has been shown that the form of the operator which intertwines the coordinate multiplication acting on holomorphic Hilbert space with matrix-valued reproducing kernels.

It is well known that unitary operators maintain rigidity, while general invertible operators destroy rigidity.
Taking this into account, we expect that the study of operator similarity is challenging, even in one variable.
The model theorem is given in Chapter 0.2 of \cite{NKN}, in the view of complex geometry, and shows that the eigenvector bundle induced by contraction in $\mathcal{B}^{1}_n(\Omega)$ has a kind of tensor structure.
By using the main result of \cite{Tr} and the model theorem for contractions, H. Kwon and S. Treil proved a theorem which allows one to decide whether  a contractive
operator $T$ is similar to the $n$ times copies of $M^*_z$ on Hardy space or not, which is  $$\Big\|\frac{\partial P(w)}{\partial w}\Big\|^2_{HS}-\frac{n}{(1-|w|^2)^2}\leq \frac{\partial^{2}}{\partial w\partial \overline{}{w}} \psi(w),
\ w\in \mathbb{D}$$
for projection-valued function $P$ with $\mbox{ran} P(w)=\mbox{ker}(T-w)$ and a bounded subharmonic function $\psi$.
Then, the result was generalized to the case of weighted Bergman shift by R.G. Douglas, H. Kwon and S. Treil \cite{Kwon1}.
Subsequently, the quantity $-\|\frac{\partial P(w)}{\partial w}\|^2_{HS}$ has been proved to be the trace of the curvature of $T$ (cf. \cite{HJK}).
Currently, this result does not have a version for commuting m-tuples.
Although there exist plenty of model theorems about the commuting operator tuples \cite{AE,AA,AA2,MV}, the techniques cannot
be easily generalized for the lack of proper condition for the Corona theorem in several variables.

In infinite-dimensional separable Hilbert spaces, strongly irreducible operators can be regarded as a natural generalization of Jordan block matrix. Strong irreducibility is a similarity invariant of operators.
In \cite{CFJ}, Y. Cao, J.S. Fang and C.L. Jiang  introduced the $K_0$-group into the similarity classification of operators and
characterized when  operators have a unique strongly irreducible decomposition up to similarity. 
Consequently, C.L. Jiang, X.Z. Guo and the second author gave a similarity theorem of Cowen-Douglas operators by using the
ordered $K$-group of the commutant algebra as an invariant \cite{JGJ}.
From the perspective of complex geometry, the similarity of Cowen-Douglas operators may be described through the equivalence of two families of eigenvectors in \cite{J}.
Using the eigenvector bundle associated to $T\in \mathcal{B}_{n}^{1}(\Omega)$, M. Uchiyama discussed when $T$ is similar or quasi-similar to the unilateral shift \cite{MU}.

In \cite{WWH}, W.W. Hastings provided a function-theoretic characterization of subnormal tuples quasi-similar to the Cauchy tuple.
Concerning absolute equivalence, virtual unitary equivalence, and almost unitarily equivalence of tuples, readers are referred to \cite{CG, AJ, EJT}.

In 2007, R.G. Douglas raised an open question \cite{Douglas} (Question 4), which has not been completely solved so far. The open question is the following. 

\textbf{Question: }Can one give conditions involving the curvatures which imply that two quasi-free Hilbert modules of multiplicity one are similar?

In this note, the main result is above the geometric similarity invariant of arbitrary Cowen-Douglas tuples without the assumptions of $n$-hypercontraction and the help of the Corona theorem. To some extent, it gives a partial answer to the question above.

The paper is organized as follows. In section 2, we recall some notions and basic results above tuples in the Cowen-Douglas class.
In section 3, we obtain  an equivalence condition for the similarity of commuting operator tuples.
Furthermore, a similarity classification theorem for tuples in $\mathcal{B}_1^m(\Omega)$ is given
by using the local equivalence of the holomorphic bundles associated with some Cowen-Douglas tuples of index two.
In section 4, we introduce a new class of commuting tuples in the Cowen-Douglas class
(notice that the unitary  intertwining operator is not diagonal in this case).
In section 5, some weakly homogeneous operators are investigated.

\section{Preliminaries}
In this section, we will recall some notations and basic results of tuples in the Cowen-Douglas class.
Let ${\mathcal L}({\mathcal H})^{m}$ be the collection of all commuting $m$-tuples of bounded operators on ${\mathcal H}$.
For $\textbf{T}=(T_{1},\cdots, T_{m})\in{\mathcal L}({\mathcal H})^{m}$, we define $\textbf{T}x=(T_1x, \cdots, T_mx)$, $x\in \H$ and $\textbf{T}-w=(T_1-w_1,\cdots, T_m-w_m)$, then
$ \mbox{ker}(\textbf{T}-w)=\bigcap_{i=1}^m\mbox{ker}(T_i-w_i)$
with $w=(w_1, \cdots, w_m)$ in $\Omega$.
The class of Cowen-Douglas tuple of operators with rank $n$ over $\Omega$: $\mathcal{B}^{m}_n(\Omega)$ is defined as follows \cite{CD,CD3}:
$$\begin{array}{lll}
\mathcal{B}^{m}_n(\Omega):=\{\textbf{T}\in{\mathcal L}({\mathcal H})^{m}:
&(1)\,\,\bigvee_{w\in \Omega}\mbox{ker}(\textbf{T}-w)=\mathcal{H},\\
&(2)\,\,\mbox{ran}(\textbf{T}-w)\ \mbox{is\ closed}\ \mbox{for\ all}~w\in\Omega,\\
&(3)\,\,\mbox{dim ker}(\textbf{T}-w)=n\ \mbox{for\ all}~w\in\Omega\}.
\end{array}$$

It follows that for each $w\in \Omega$, $\ker(\textbf{T}-w)$ is an $n$-dimensional vector subspace of $\mathcal{H}$. Define $E_\textbf{T}:=\{(w,x)\in\Omega\times\mathcal{H}:x\in\mbox{ker}(\textbf{T}-w)\}$ with a projection map $\pi:E_\textbf{T}\rightarrow\Omega$ such that $\pi^{-1} (w)=\ker(\textbf{T}-w)$.
It is a sub-bundle of $\Omega\times\mathcal{H}$ and its Hermitian structure comes from $\mathcal{H}$.
Thus, $E_{\textbf{T}}$ associated with $\textbf{T}$ is an $n$-dimensional Hermitian holomorphic vector bundle.
\begin{thm}\cite{CD,CD3}\label{2110.11}
Let $\text{\bf{T}}, \text{\bf{S}}\in \mathcal{B}^{m}_n(\Omega)$. Then $\text{\bf{T}}\sim_{u}\text{\bf{S}}$ if and only if the Hermitian holomorphic vector bundles
$E_{\text{\bf{T}}}$ and $E_{\text{\bf{S}}}$ are congruent $( \mbox{denoted\, by}\  E_\text{\bf{T}}\sim_{u}
E_\text{\bf{S}})$ over some open subset $\Omega_0$ of $\Omega\subset \mathbb{C}^{m}$.
\end{thm}

When $m=1$, the above theorem is proved in Theorem 2.6 of \cite{CD}, in the case of $m>1$, Theorem \ref{2110.11} is also valid (\cite{CD3}, pp.16) due to M.J. Cowen and R.G. Douglas.
They make a rather detailed study
of certain aspects of complex geometry and introduce the following concepts.

Let $E$ be a $C^\infty$ vector bundle over $\Omega$. A connection $D$ is a differential operator, which takes sections of $E$ to sections with 1-form coefficients and satisfies the Leibnitz rule $D(fs)=(df)s+fDs$ for section $s$ and function $f$.
Similarly, $D^2$ can be defined, $D^2s=\mathcal{K}s dzd\bar{z}$ for section $s$, bundle map $\mathcal{K}$ determined by $D^2$ is called as the curvature of bundle $E$.

For every Hermitian holomorphic vector bundle $E$ over $\Omega$, there is a unique canonical connection $\Theta$, which is a Chern connection metric-preserving and compatible with the holomorphic structure. 
Given a holomorphic frame $\gamma=\{\gamma_i\}_{i=1}^{n}$ of $E$, we have the metric
$h(w)=(\!(\langle\gamma_j(w),\gamma_i(w)\rangle)\!)_{n\times n}$ and $D\gamma=\gamma\Theta$, $\Theta=(\Theta_{ij})_{i,j=1}^{n}$ is the matrix of connection 1-form.
The curvature of $E$ can be defined as:
\begin{equation}\label{2110.111}
\mathcal{K}(w)
=-\sum_{i,j=1}^m\frac{\partial}{\partial \bar{w}_{j}}\Big(h^{-1}(w)\frac{\partial
h(w)}{\partial w_{i}}\Big)dw_i\wedge d\bar{w}_{j}
\end{equation}
for $w=(w_1,\cdots,w_m)\in\Omega$.
When $E$ is a line bundle, equation (\ref{2110.111}) is equivalent to $\mathcal{K}(w)
=-\sum_{i,j=1}^m\frac{\partial^{2}\log\|\gamma(w)\|^2}{\partial \bar{w}_{j}\partial w_i}dw_i\wedge d\bar{w}_{j}$, where $\gamma$ is a non-zero section of $E$.

For any $C^{\infty}$ bundle map $\phi$ on $E$ and given frame $\sigma$ of $E$, we have that

$\begin{array}{llll}
 &(1)\,\,\phi_{\bar{w}}(\sigma)=\frac{\partial}{\partial \bar{w}}(\phi(\sigma));\\
 &(2)\,\,\phi_{w}(\sigma)=\frac{\partial}{\partial w}(\phi(\sigma))+[h^{-1}\frac{\partial}{\partial w}h,\phi(\sigma)].
\end{array}$

Since the curvature can also be regarded as a bundle map, we obtain covariant derivatives ${\mathcal K}_{w^I\bar{w}^J}$, $I,J\in\mathbf{Z}_{+}^{m}$ of the curvature by using the
inductive formulaes above, where $\mathbf{Z}_{+}^{m}$ is the collection of $m$-tuples of nonnegative integers.
The curvature ${\mathcal K}$ and it's covariant derivatives ${\mathcal K}_{w^I\bar{w}^J}$ are the
unitarily invariants of Hermitian holomorphic vector bundle $E$ (see \cite{CD,CD3}).

\begin{thm}\cite{CD,CD3}\label{3.5you}
Let $\text{\bf{T}}, \text{\bf{S}}\in \mathcal{B}^{m}_n(\Omega)$. Then $E_\text{\bf{T}}\sim_{u}
E_\text{\bf{S}}$ if and only if there exists an isometry $V: E_{\text{\bf{T}}}\rightarrow E_{\text{\bf{S}}}$ and a number $m$ depending on $E_{\text{\bf{T}}},E_{\text{\bf{S}}}$ such that
$$V{\mathcal K}_{\text{\bf{T}},w^I\bar{w}^J}={\mathcal K}_{\text{\bf{S}},w^I\bar{w}^J}V, I,J\in\mathbf{Z}_{+}^{m}.$$
\end{thm}

\section{On the similarity of commuting operator tuples}
The classification of similarities of commuting operator tuples has always been a challenging problem.
Even in the operator case, it is not yet clear how to describe the similarity of Cowen-Douglas operators in $\mathcal{B}^{1}_1(\Omega)$ using only geometric quantities, such as the curvature. M.J. Cowen and R.G. Douglas put forward the
following conjecture in 4.35 of \cite{CD}:
If $\bar{\mathbb{D}}$ (the closure of unit disc $\mathbb{D}$) is a $k$-spectral set for $T,S\in \mathcal{B}_1^1(\mathbb{D})$, then $T\sim_{s}S$ if and only if their curvatures $\mathcal{K}_{T}$ and $\mathcal{K}_{S}$ satisfy
$\lim\limits_{w\rightarrow \partial\mathbb{D}}\frac{\mathcal{K}_{T}(w)}{\mathcal{K}_{S}(w)}=1$.
In \cite{CM1,CM2}, two counter examples were constructed by D.N. Clark
and G. Misra.
Instead of the quotient of the curvatures, they considered the
quotient of metrics $h_T$ and $h_S$ of $E_T$ and $E_S$ denoted by
$a_{w}$.  It was then proved in \cite{CM2} that contraction
$T$ is similar to $S_{\alpha}$ (with weight sequence $\{(\frac{n+1}{n+2})^{\frac{\alpha}{2}}\}_{n=0}^{\infty}$) if and only if $a_{w}$ is bounded and bounded below by $0$.
This result can be regarded as a geometric version of the classical result for the weighted shifts given by A.L. Shields (see \cite{SAL}).
For recent developments concerning the similarity of Cowen-Douglas operators, the reader is referred to \cite{DHS,DKKS1, DKKS2,JJKMCR,Kwon1}.

Although there are many studies on the similarity classification of Cowen-Douglas operators, the similarity classification of commuting tuples  is not yet fully solved.  In this chapter, we provide a different necessary and sufficient condition for the similarity of commuting operator tuples.
We introduce the following definition of $\sigma_{\textbf{T}_0,\textbf{T}_1}$, and the notation is adopted from the next.

\begin{defn} Let $\textbf{T}_i\in{\mathcal L}({\mathcal H}_i)^{m},i=0,1$. Define $\sigma_{\textbf{T}_0,\textbf{T}_1}:{\mathcal L}({\mathcal H}_1,{\mathcal H}_0)\rightarrow
{\mathcal L}({\mathcal H}_{0})^{m}$ be the tuple
$$\sigma_{\textbf{T}_0,\textbf{T}_1}(X)=\textbf{T}_0X-X\textbf{T}_1,\,\, X\in {\mathcal L}({\mathcal H}_1,{\mathcal H}_0).$$ Let $\sigma_{\textbf{T}_0}:{\mathcal L}({\mathcal
H}_{0})\rightarrow {\mathcal L}({\mathcal H}_0)^{m}$  be the tuple $\sigma_{\textbf{T}_0, \textbf{T}_0}.$
\end{defn}
\subsection{On the similarity of commuting tuples}

In order to describe clearly Theorem \ref{main9.26}, Lemma \ref{mainlemma}, Corollary \ref{2110.26} and Theorem \ref{main2},
we need to introduce the following notations.
Unless otherwise specified, we always assume that
$$\textbf{T}_{ij}=(T_{ij}^1,\cdots, T_{ij}^m),\ \widetilde{\textbf{T}}_{ij}=(\widetilde{T}_{ij}^1, \cdots,\widetilde{T}_{ij}^m),\ \textbf{S}_{ij}=(S_{ij}^1, \cdots, S_{ij}^m), 0\leq i\leq j\leq1$$ and
$$\textbf{T}=(T_{11}^1, \cdots, T_{11}^m),\ \textbf{S}=(S_{00}^1, \cdots, S_{00}^m)$$
for some positive integer $m$. The main theorem of this paper as follows.

\begin{thm}\label{main9.26}
Let $\text{\bf{T}},\ \text{\bf{S}}\in \mathcal{L}(\mathcal{H})^{m}$.
Suppose that $\{\text{\bf{S}}_{11}\in \mathcal{L}(\mathcal{H})^m: \mbox{ker}\ \sigma_{\text{\bf{S}}_{11},\text{\bf{T}}}=\{0\}\}\neq \emptyset$.
Then $\text{\bf{T}}\sim_{s}\text{\bf{S}}$ if and only if there exist two operator tuples $\widetilde{\text{\bf{T}}}=(T_{1}, \cdots, T_{m}), \widetilde{\text{\bf{S}}}=(S_{1}, \cdots, S_{m})\in \mathcal{L}(\mathcal{H}\oplus \mathcal{H})^{m}$ such that
\begin{enumerate}
  \item [(1)] $T_{i}=\left (\begin{smallmatrix}
 T_{00}^i & T_{01}^i\\
 0 & T_{11}^i\\
\end{smallmatrix}\right ),S_{i}=\left(\begin{smallmatrix}
 S_{00}^i & S_{01}^i \\
 0 & S_{11}^i \\
\end{smallmatrix}\right),1\leq i\leq m, $ where $\text{\bf{T}}_{01}\in \mbox{ran}\,\sigma_{\text{\bf{T}}_{00},\text{\bf{T}}},\text{\bf{S}}_{01}\in \mbox{ran}\,\sigma_{\text{\bf{S}},\text{\bf{S}}_{11}}$ and $\mbox{ker}\ \sigma_{\text{\bf{T}}_{00},\text{\bf{S}}}=\{0\}$;
  \item [(2)] $\widetilde{\text{\bf{T}}}\sim_{u}\widetilde{\text{\bf{S}}}$. 
\end{enumerate}
\end{thm}

In order to prove our main theorem, we first need a lemma which characterizes the unitary operator which intertwines two special commuting operator tuples.

\begin{lem} \label{mainlemma}
Let $\widetilde{\text{\bf{T}}}=(T_{1}, \cdots, T_{m}), \widetilde{\text{\bf{S}}}=(S_{1}, \cdots, S_{m})\in \mathcal{L}(\mathcal{H}\oplus \mathcal{H})^{m}$, where $T_{i}=\left (\begin{smallmatrix}
 T_{00}^i & T_{01}^i\\
 0 & T_{11}^i\\
\end{smallmatrix}\right ),S_{i}=\left(\begin{smallmatrix}
 S_{00}^i & S_{01}^i \\
 0 & S_{11}^i \\
\end{smallmatrix}\right),1\leq i\leq m,$ and $\text{\bf{T}}_{01}= -\sigma_{\text{\bf{T}}_{00},\text{\bf{T}}}(X),\text{\bf{S}}_{01}=-\sigma_{\text{\bf{S}},\text{\bf{S}}_{11}}(Y)$ for some $X,Y\in\mathcal{L}(\mathcal{H})$.
Suppose that $\mbox{ker}\ \sigma_{\text{\bf{T}}_{00},\text{\bf{S}}}=\mbox{ker}\ \sigma_{\text{\bf{S}}_{11},\text{\bf{T}}}=\{0\}$,
then there exists a unitary operator $U=(\!(U_{i,j})\!)_{2\times 2}$ such that $U\widetilde{\text{\bf{T}}}=\widetilde{\text{\bf{S}}}U$ if and only if
the following statements hold
 \begin{enumerate}
\item $U_{10}T_{00}^iU^{-1}_{10}=S_{11}^i,\,U^{*-1}_{01}T_{11}^iU^*_{01}=S_{00}^i,\,1\leq i\leq m$;
\item $(I+XX^*)^{-1}=U^*_{10}U_{10}, (I+X^*X)^{-1}=U^*_{01}U_{01}$;
\item $Y-U_{01}X^*U^{-1}_{10}\in ker\ \sigma_{\text{\bf{S}},\text{\bf{S}}_{11}}.$
\end{enumerate}
\end{lem}

\begin{proof}
Let $U=\left(\begin{smallmatrix}
U_{00} & U_{01}\\
U_{10} & U_{11} \\
\end{smallmatrix}\right)$. From $U\widetilde{\textbf{T}}=\widetilde{\textbf{S}}U$, we have
\begin{equation}\label{mp1}
U_{10}XT_{11}^i-U_{10}T_{00}^iX=S_{11}^iU_{11}-U_{11}T_{11}^i,
\end{equation}
\begin{equation}\label{mp2}
U_{00}T_{00}^i-YS_{11}^iU_{10}=S_{00}^iU_{00}-S_{00}^iYU_{10}, T_{00}^iU^*_{00}+(XT_{11}^i-T_{00}^iX)U^*_{01}=U^*_{00}S_{00}^i
\end{equation}
and
\begin{eqnarray}\label{mp3}
U_{10}T_{00}^i=S_{11}^iU_{10}, T_{11}^iU^*_{01}=U^*_{01}S_{00}^i,\,1\leq i\leq m.
\end{eqnarray}

First of all, we will prove that $U_{01}$ and $U_{10}$ are invertible.
By (\ref{mp1}) and (\ref{mp3}), we have
$U_{10}XT_{11}^i-S_{11}^iU_{10}X=S_{11}^iU_{11}-U_{11}T_{11}^i$ and
$(U_{10}X+U_{11})T_{11}^i=S_{11}^i(U_{10}X+U_{11}),\,1\leq i\leq m.$
From (\ref{mp2}) and (\ref{mp3}), we also have
$T_{00}^i(U^*_{00}-XU^*_{01})=(U^*_{00}-XU^*_{01})S_{00}^i,\,1\leq i\leq m.$
It follows that
$$U_{10}X+U_{11}\in \cap_{i=1}^{m}\mbox{ker}\ \sigma_{S_{11}^i,T_{11}^i}=\mbox{ker}\ \sigma_{\textbf{S}_{11},\textbf{T}},\,U^*_{00}-XU^*_{01}\in \cap_{i=1}^{m}\mbox{ker}\ \sigma_{T_{00}^i,S_{00}^i}=\mbox{ker}\ \sigma_{\textbf{T}_{00},\textbf{S}}.$$
Note that $\mbox{ker}\ \sigma_{\textbf{S}_{11},\textbf{T}}=\mbox{ker}\ \sigma_{\textbf{T}_{00},\textbf{S}}=\{0\}$. We see
\begin{eqnarray}\label{mp0614}U_{00}=U_{01}X^*, U_{11}=-U_{10}X.\end{eqnarray}
So the form of the unitary operator $U$ is
$\left (\begin{smallmatrix} U_{01}X^*&U_{01}\\U_{10}&-U_{10}X \end{smallmatrix}\right).$
By using the fact $UU^*=U^*U=I\oplus I$, we have the following equations:
\begin{eqnarray}\label{mp7}U_{01}(I+X^*X)U^*_{01}=I,U_{10}(I+XX^*)U^*_{10}=I,\end{eqnarray}
\begin{eqnarray}\label{mp9}XU^*_{01}U_{01}=U^*_{10}U_{10}X,\end{eqnarray}
\begin{eqnarray}\label{mp10}XU^*_{01}U_{01}X^*+U^*_{10}U_{10}=I, X^*U^*_{10}U_{10}X+U^*_{01}U_{01}=I.\end{eqnarray}
By equations (\ref{mp7})-(\ref{mp10}), we have $U^*_{10}U_{10}(I+XX^*)=I$ and $(I+X^*X)U^*_{01}U_{01}=I$.
Moreover, combining equation (\ref{mp7}), we obtain $U_{01}$ and $U^*_{10}$ are invertible and also $(I+X^*X)U^*_{01}$ and $U_{10}(I+XX^*)$.
Since $I+XX^*$ and $I+X^*X$ are invertible, it is easy to see that
$U^*_{10}U_{10}=(I+XX^*)^{-1}$ and $U^*_{01}U_{01}=(I+X^*X)^{-1}$. From equation (\ref{mp3}) and the invertibility of $U_{01}$ and $U_{10}$, we imply  $U_{10}T_{00}^iU^{-1}_{10}=S_{11}^i$ and $U^{*-1}_{01}T_{11}^iU^*_{01}=S_{00}^i$ for all $1\leq i\leq m$.

By equation (\ref{mp2}), for any $1\leq i\leq m$, we have
\begin{eqnarray}\label{mp12}
S_{00}^iYU_{10}-S_{00}^iU_{00}=YS_{11}^iU_{10}-U_{00}T_{00}^i=YS_{11}^iU_{10}-U_{00}U^{-1}_{10}S_{11}^iU_{10}.
\end{eqnarray}
Multiplying $U^{-1}_{10}$ on the right side of equation (\ref{mp12}), we obtain
$S_{00}^iY-S_{00}^iU_{00}U^{-1}_{10}=YS_{11}^i-U_{00}U^{-1}_{10}S_{11}^i,\,1\leq i\leq m.$
Since equation (\ref{mp0614}), it follows that
$S_{00}^i(Y-U_{01}X^*U_{10}^{-1})=(Y-U_{01}X^*U_{10}^{-1})S_{11}^i,\,1\leq i\leq m.$
That is, $Y-U_{01}X^*U^{-1}_{10}\in \cap_{i=1}^{m}\mbox{ker}\ \sigma_{S_{00}^i,S_{11}^i}=ker\ \sigma_{\textbf{S},\textbf{S}_{11}}.$

For the sufficient part, let $U=\left (\begin{smallmatrix} U_{01}X^*&U_{01}\\U_{10}&-U_{10}X \end{smallmatrix}\right)$ which satisfies the conditions (1)-(3). It implies $U$ is a unitary operator.
In the following, we will check $U\widetilde{\textbf{T}}U^*=\widetilde{\textbf{S}}$, that is, $UT_{i}U^*=S_{i},\,1\leq i\leq m.$
Note that $UT_{i}U^*$ has the following form
$$
\left (\begin{smallmatrix}U_{01}X^*T_{00}^iXU^*_{01}+\big(U_{01}X^*(XT_{11}^i-T_{00}^iX)+U_{01}T_{11}^i\big)U_{01}^{*} &U_{01}X^*T_{00}^iU^*_{10}+\big(U_{01}X^*(XT_{11}^i-T_{00}^iX)+U_{01}T_{11}^i\big)(-X^*U^*_{10})\\
U_{10}T_{00}^iXU^*_{01}+U_{10}(XT_{11}^i-T_{00}^iX)U^*_{01}-U_{10}XT_{11}^iU^*_{01}&U_{10}T_{00}^iU^*_{10}+U_{10}(XT_{11}^i-T_{00}^iX)(-X^*U^*_{10})+U_{10}XT_{11}^iX^*U^*_{10}\\
\end{smallmatrix}\right).$$
Since $Y-U_{01}X^*U^{-1}_{10}\in \mbox{ker}\ \sigma_{\textbf{S},\textbf{S}_{11}},$
$YS_{11}^i-S_{00}^iY=(U_{01}X^*U^{-1}_{10})S_{11}^i-S_{00}^i(U_{01}X^*U^{-1}_{10}), 1\leq i\leq m.$
By statements (1) and (2), we have
\begin{eqnarray*}
   &&U_{01}X^*T_{00}^iU^*_{10}+\big(U_{01}X^*(XT_{11}^i-T_{00}^iX)+U_{01}T_{11}^i\big)(-X^*U^*_{10})\\
   &=&U_{01}X^*T_{00}^i(I+XX^*)U_{10}^*-U_{01}(I+X^*X)T_{11}^iX^*U^*_{10}\\
    &=&U_{01}X^*U_{10}^{-1}S_{11}^i-S_{00}^iU^{*-1}_{01}X^*U^*_{10}\\
   &=&YS_{11}^i-S_{00}^iY,\,1\leq i\leq m.
\end{eqnarray*}
Based on a routine computation, we obtain $U\widetilde{\textbf{T}}U^*=\widetilde{\textbf{S}}$.
The proof of these equations also use that fact
$X^*U^*_{10}U_{10}=X^*(I+XX^*)^{-1}=(I+X^*X)^{-1}X^*=U^*_{01}U_{01}X^*.$
These equalities finish the proof of sufficient part.
\end{proof}
\textbf{The proof of Theorem \ref{main9.26}: }
\begin{proof}
Sufficiency:
Firstly, there is a tuple $\textbf{S}_{11}$ such that $\mbox{ker}\ \sigma_{\textbf{S}_{11},\textbf{T}}=\{0\}$.
Suppose that two commuting tuples of operators $\widetilde{\textbf{T}}=(T_{1}, \cdots, T_{m})$ and $\widetilde{\textbf{S}}=(S_{1}, \cdots, S_{m})\in \mathcal{L}(\mathcal{H}\oplus \mathcal{H})^{m}$ are unitarily equivalent, that is, there exists a
unitary operator $U=\left(\begin{smallmatrix}
U_{00} & U_{01}\\
U_{10} & U_{11} \\
\end{smallmatrix}\right)$ such that $U\widetilde{\textbf{T}}=\widetilde{\textbf{S}}U$.
If condition (1) in the theorem is satisfied by $\widetilde{\textbf{T}},\widetilde{\textbf{S}}$,
by Lemma \ref{mainlemma}, we know that $U_{01}^{*}$ is invertible and $T_{11}^iU_{01}^{*}=U_{01}^{*}S_{00}^i$ for all $1\leq i\leq m$.
It follows that $\textbf{T}U_{01}^{*}=U_{01}^{*}\textbf{S}$.
Thus, $\textbf{T}$ is similar to $\textbf{S}$.

Necessity: Since $\textbf{T}$ is similar to $\textbf{S}$, there exists an invertible operator $X_{1}$ such that
\begin{equation}\label{3.1011}
\textbf{T}=X_{1}\textbf{S}X_{1}^{-1},
\end{equation}
that is, $T_{11}^i=X_{1}S_{00}^iX_{1}^{-1}$ for all $1\leq i\leq m$.
Without loss of generality, we assume that $(X_{1}^{-1})^{*}X_{1}^{-1}-I\geq0$. Otherwise, let $\alpha=\inf\{x|x\in \sigma\big((X_{1}^{-1})^{*}X_{1}^{-1}\big)\}$,
we have $\frac{(X_{1}^{-1})^{*}X_{1}^{-1}}{\alpha}-I\geq0,$ since $X_{1}$ is invertible.
Then notice that $\alpha>0$, upon replacing $X_{1}^{-1}$ by $\frac{X_{1}^{-1}}{\sqrt{\alpha}}$, we obtain that $(X_{1}^{-1})^{*}X_{1}^{-1}-I\geq0$.
Therefore, we find a bounded linear operator $X$, such that
\begin{equation}\label{3.1021}
I+X^{*}X=(X_{1}^{-1})^{*}X_{1}^{-1}.
\end{equation}
Obviously, $I+XX^{*}$ is also invertible and positive. In the same way as constructing $X$, we know that there exists $X_2$ satisfies
\begin{equation}\label{3.1031}
(I+XX^{*})^{-1}=X_2^{*}X_2,
\end{equation}
and $X_2$ is an invertible operator.

Choosing a non-zero commuting tuple of operators $\textbf{S}_{11}\in \mathcal{L}(\mathcal{H})^{m}$ such that  $\mbox{ker}\ \sigma_{\textbf{S}_{11},\textbf{T}}=\{0\}$, that is,
$\bigcap_{i=1}^{m}\mbox{ker}\ \sigma_{S_{11}^i,T_{11}^i}=\{0\}$.
Next, we will construct another tuple $\textbf{T}_{00}\in \mathcal{L}(\mathcal{H})^{m}$.
Let
\begin{equation}\label{3.1041}
T_{00}^i:=X^{-1}_2S_{11}^iX_2,
\end{equation}
for $1\leq i\leq m$. Then $\textbf{T}_{00}=X^{-1}_2\textbf{S}_{11}X_2$.

We claim that $\mbox{ker}\ \sigma_{\textbf{T}_{00},\textbf{S}}=\{0\}$.
If $Z\in \mbox{ker}\ \sigma_{\textbf{T}_{00},\textbf{S}}$, then for all $1\leq i\leq m$, we have $T_{00}^iZ=ZS_{00}^i$, equivalently,
$S_{11}^iX_{2}ZX_{1}^{-1}=X_{2}ZX_{1}^{-1}T_{11}^i$, since equations (\ref{3.1011}) and (\ref{3.1041}) hold.
By $\bigcap_{i=1}^{m}\mbox{ker}\ \sigma_{S_{11}^i,T_{11}^i}=\{0\}$, we obtain $X_{2}ZX_{1}^{-1}=0$.
Note that $X_1, X_2$ are both invertible. It follows that $Z=0$ and $\bigcap_{i=1}^{m}\mbox{ker}\ \sigma_{T_{00}^i,S_{00}^i}=\mbox{ker}\ \sigma_{\textbf{T}_{00},\textbf{S}}=\{0\}.$

For any bounded linear operator $W\in \bigcap_{i=1}^{m}\mbox{ker}\ \sigma_{S_{00}^i,S_{11}^i}$, let $Y:=W+X_1^{*}X^{*}X_2^{-1}$. This implies that
\begin{equation}\label{3.1051}
Y-X_1^{*}X^{*}X_2^{-1}\in \cap_{i=1}^{m}\mbox{ker}\ \sigma_{S_{00}^i,S_{11}^i}=\mbox{ker}\ \sigma_{\textbf{S},\textbf{S}_{11}}.
\end{equation}
Based on the above discussion, we may assume that operator tuples $\widetilde{\textbf{T}}=(T_{1}, \cdots, T_{m}), \widetilde{\textbf{S}}=(S_{1}, \cdots, S_{m})$ with
$T_{i}=\left (\begin{smallmatrix}
 T_{00}^i & T_{01}^i\\
 0 & T_{11}^i\\
\end{smallmatrix}\right ),S_{i}=\left(\begin{smallmatrix}
 S_{00}^i & S_{01}^i \\
 0 & S_{11}^i \\
\end{smallmatrix}\right),1\leq i\leq m$ and $\textbf{T}_{01}= -\sigma_{\textbf{T}_{00},\textbf{T}}(X),\textbf{S}_{01}=-\sigma_{\textbf{S},\textbf{S}_{11}}(Y)$.
Then a simple calculation shows that $\widetilde{\textbf{T}}$ and $\widetilde{\textbf{S}}$ are commuting tuples and satisfy the condition (1).

Set $U:=\left(\begin{smallmatrix}X_{1}^*X^*&X_{1}^*\\X_{2}&-X_{2}X\end{smallmatrix}\right)$. From Lemma \ref{mainlemma} and equations (\ref{3.1011})-(\ref{3.1051}), we obtain $U$ is unitary and $U\widetilde{\textbf{T}}=\widetilde{\textbf{S}}U$. Hence, $\widetilde{\textbf{T}}\sim_{u}\widetilde{\textbf{S}}$.
\end{proof}

Given an $m$-tuple of operators $\mathbf{T}=(T_1,\cdots,T_m)\in\mathcal{B}^{m}_n(\Omega)$, by subsection 2.2 in \cite{Dinesh}, we know that $\mathbf{T}$ is unitarily equivalent to the adjoint of an $m$-tuple of multiplication operators $\mathbf{M}_z=(M_{z_1},\cdots,M_{z_m})$ by coordinate functions on some Hilbert
space $\mathcal{H}_K$ of holomorphic functions on $\Omega^*=\{\bar{w}:w\in\Omega\}$ possessing
a reproducing kernel $K$. It is expressed equivalently as $\mathbf{T}\sim_{u}(\mathbf{M}_z^*,\mathcal{H}_K)$.
Define $e_w$ be the evaluation function of $\mathcal{H}_K$ at $w$. Given a vector $\xi\in \mathbb{C}^n$, the function $e^*_{\bar{w}}\xi\in\mathcal{H}_K$
and is denoted by $K(\cdot, \bar{w})\xi$, which has the reproducing property $\langle f,K(\cdot,\bar{w})\xi\rangle_{\mathcal{H}_K}=\langle f(\bar{w}),\xi\rangle_{\mathbb{C}^n}$.
In addition, we have $\ker(\mathbf{M}_z^*-w)=\{K(\cdot,\bar w)\xi,\xi\in \mathbb{C}^n\}$.

In order to find the minimal order $m$ of covariant partial derivatives in Theorem \ref{3.5you}, M.J. Cowen and R.G. Douglas introduced the concept of coalescing set \cite{CD2}.
The algebra $\mathscr{A}(w)$ is generated by the curvatures and their covariant derivatives at $w$.
The coalescing set of $\mathscr{A}(x)$ is the set where the dimension of $\mathscr{A}(x)$ (as a function of $x$) is not locally constant.
It is trivially closed and nowhere dense.
Furthermore, they proved that two bundles are locally equivalent on at least one dense open set, i.e. the complement of the coalescing set for the curvature corresponding to one of bundles,
 which means that the two bundles are equivalent.
In Corollary \ref{2110.26} and Theorem \ref{main2}, using this geometric quantity,
we characterize the similarity classification of tuples in the Cowen-Douglas class.
In other words, we give a partial answer to R.G. Douglas's question about the geometric similarity of Cowen-Douglas class for multivariable case.
Our results allow one to use geometric quantities of Cowen-Douglas tuples with index one to assess whether they are similar or not.

\begin{cor}\label{2110.26}
Let $\text{\bf{T}}, \text{\bf{S}}\in \mathcal{B}^{m}_n(\Omega)$
and $\text{\bf{T}}\sim_{u} (\text{\bf{M}}^*_z, \mathcal{H}_{K_{\text{\bf{T}}}}),\text{\bf{S}}\sim_{u} (\text{\bf{M}}^*_z, \mathcal{H}_{K_{\text{\bf{S}}}}).$
Suppose that $\{\text{\bf{S}}_{11}\in \mathcal{L}(\mathcal{H})^m: \mbox{ker}\ \sigma_{\text{\bf{S}}_{11},\text{\bf{T}}}=\{0\}\}\neq \emptyset$.
Then $\text{\bf{T}}\sim_{s}\text{\bf{S}}$ if and only if there exist two tuples $\widetilde{\text{\bf{T}}}=(T_{1}, \cdots, T_{m}), \widetilde{\text{\bf{S}}}=(S_{1}, \cdots, S_{m})\in \mathcal{B}^{m}_{2n}(\Omega)$ such that
\begin{enumerate}
  \item [(1)]
$T_{i}=\left (\begin{smallmatrix}
 T_{00}^i & T_{01}^i\\
 0 & T_{11}^i\\
\end{smallmatrix}\right ),S_{i}=\left(\begin{smallmatrix}
 S_{00}^i & S_{01}^i \\
 0 & S_{11}^i \\
\end{smallmatrix}\right),1\leq i\leq m$, where
$\text{\bf{T}}_{01}\in \mbox{ran}\,\sigma_{\text{\bf{T}}_{00},\text{\bf{T}}},\text{\bf{S}}_{01}\in \mbox{ran}\,\sigma_{\text{\bf{S}},\text{\bf{S}}_{11}}$ and $\mbox{ker}\ \sigma_{\text{\bf{T}}_{00},\text{\bf{S}}}=\{0\}$;
  \item [(2)] the bundles $E_{\widetilde{\text{\bf{T}}}}$ and $E_{\widetilde{\text{\bf{S}}}}$ of $\widetilde{\text{\bf{T}}}$ and $\widetilde{\text{\bf{S}}}$ are locally equivalent on an open dense subset of $\Omega$, the complement of the coalescing set for the curvature of $E_{\widetilde{\text{\bf{T}}}}$.  
\end{enumerate}
\end{cor}

\begin{proof}
Let $\Omega_{0}$ be the complement of the coalescing set for the curvature of $E_{\widetilde{\textbf{T}}}$. Clearly, $\Omega_{0}\subset\Omega$.
If the bundles $E_{\widetilde{\textbf{T}}}$ and $E_{\widetilde{\textbf{S}}}$ are locally equivalent on $\Omega_{0}$,
by using the main theorem of \cite{CD2} due to M.J. Cowen and R.G. Douglas,
we obtain that
  metric-preserving connections $D_{\widetilde{\textbf{T}}}$ and $D_{\widetilde{\textbf{S}}}$ of $E_{\widetilde{\textbf{T}}}$ and $E_{\widetilde{\textbf{S}}}$ are equivalent to order $2n$ on $\Omega$.
From the proof of the sufficiency in Theorem \ref{main9.26}, we see that
$\textbf{T}$ is similar to $\textbf{S}$.

From the proof of Theorem \ref{main9.26} and the main theorem of \cite{CD2}, we only need to prove $\widetilde{\textbf{T}}, \widetilde{\textbf{S}}\in \mathcal{B}^{m}_{2n}(\Omega)$.

Suppose that there exist $X,Y$ such that $\textbf{T}_{01}=\sigma_{\textbf{T}_{00},\textbf{T}}(-X),\textbf{S}_{01}=\sigma_{\textbf{S},\textbf{S}_{11}}(-Y)$, that is,
$T_{01}^j=XT_{11}^j-T_{00}^jX,S_{01}^j=YS_{11}^j-S_{00}^jY$ for all $1\leq j\leq m$.
Without losing generality, we assume that $\textbf{T}=(\textbf{M}^*_z, \mathcal{H}_{K_{\textbf{T}}})$, and then
$\ker(\textbf{T}-w)=\{K_{\textbf{T}}(\cdot,\bar{w})\xi,\xi\in \mathbb{C}^n\},w\in\Omega$. For fixed but arbitrary $w\in\Omega$ and $\xi\in \mathbb{C}^n$,
we have
$$\begin{array}{lll}
T_{01}^jK_{\textbf{T}}(\cdot,\bar{w})\xi
&=&(XT_{11}^j-T_{00}^jX)K_{\textbf{T}}(\cdot,\bar{w})\xi\\
&=&w_{j}XK_{\textbf{T}}(\cdot,\bar{w})\xi-T_{00}^jXK_{\textbf{T}}(\cdot,\bar{w})\xi \\
&=&(T_{00}^j-w_{j})(-XK_{\textbf{T}}(\cdot,\bar{w})\xi),1\leq j\leq m.
\end{array}$$
It follows that $T_{01}^j( \mbox{ker}(\textbf{T}-w))\subset \mbox{ran}(\textbf{T}_{00}-w),1\leq j\leq m$
and $\textbf{T}_{01}( \mbox{ker}(\textbf{T}-w))\subset \mbox{ran}(\textbf{T}_{00}-w).$
Thus, $\widetilde{\textbf{T}}\in \mathcal{B}^{m}_{2n}(\Omega)$. Similarly, we also have $\widetilde{\textbf{S}}\in \mathcal{B}^{m}_{2n}(\Omega)$.
This completes the proof.
\end{proof}


C.L. Jiang, D.K. Keshari, G. Misra and the second author in \cite{JJKM} showed that for $T,\tilde{T}\in \mathcal{B}^{1}_1(\Omega)$, if $XT=\tilde{T}X$, then either $X=0$ or $X$ has a dense range.
In fact, we see that this result is also true when $\textbf{T},\widetilde{\textbf{T}}$ are tuples in the Cowen-Douglas class with index one. The next lemma  shows that the conditions in Lemma \ref{mainlemma} can be satisfied in many cases.
\begin{lem}\label{intertwining}
Let $\text{\bf{T}},\text{\bf{S}}\in \mathcal{B}^{m}_1(\Omega)$ and $\text{\bf{T}}\sim_{u} (\text{\bf{M}}^*_z, \mathcal{H}_{K_{0}}),\text{\bf{S}}\sim_{u} (\text{\bf{M}}^*_z, \mathcal{H}_{K_{1}}).$
If $$\lim\limits_{dist(w, \partial \Omega)\rightarrow 0}\frac{K_0(w,w)}{K_1(w,w)}=0,$$ then there exists no non-zero bounded intertwining operator $X$ such that $X\text{\bf{T}}=\text{\bf{S}}X$, i.e. $ker\ \sigma_{\text{\bf{S}},\text{\bf{T}}}=\{0\}$.
\end{lem}

\begin{proof}
Without loss of generality, we set $\textbf{T}=(\textbf{M}^*_z, \mathcal{H}_{K_{0}}),\textbf{S}=(\textbf{M}^*_z, \mathcal{H}_{K_{1}}),$ where $\mathcal{H}_{K_{i}}$ are vector-valued analytic functional Hilbert spaces with reproducing kernels $K_i, i=0,1$, respectively.
Suppose that $X\textbf{T}=\textbf{S}X$ for a bounded operator $X$. This means that $X(\mbox{ker}(\textbf{T}-w))\subset\mbox{ker}(\textbf{S}-w),\,w\in\Omega.$
Since $K_0(\cdot,\bar{w})\in \mbox{ker}(\textbf{T}-w),\,K_1(\cdot,\bar{w})\in \mbox{ker}(\textbf{S}-w)$,
there exists a holomorphic function $\phi$ on $\Omega$ such that $X(K_0(\cdot,\bar{w}))=\phi(w)K_1(\cdot,\bar{w}), w\in \Omega$ (see details in Proposition 2.4 \cite{R}). Note that
$\big\|X(\frac{K_0(\cdot,\bar{w})}{\|K_0(\cdot,\bar{w})\|})\big\|=|\phi(w)|\frac{\|K_1(\cdot,\bar{w})\|}{\|K_0(\cdot,\bar{w})\|}\leq \|X\|,$
we have $|\phi(w)|\leq \|X\|\frac{\|K_0(\cdot,\bar{w})\|}{\|K_1(\cdot,\bar{w})\|}$.
By using $\lim\limits_{dist(w, \partial \Omega)\rightarrow 0}\frac{K_0(w,w)}{K_1(w,w)}=0$, we obtain that $|\phi|$ will goes to zero when $dist(w, \partial \Omega)$ goes to zero.
By the maximum modulus principle of holomorphic function, we have $\phi(w)$ is equal to zero for all $w\in \Omega$, so does $X(K_0(\cdot,\bar{w}))$. According to the spanning property $\ker(\textbf{T}-w)=\bigvee\{K_0(\cdot,\bar{w})\}$, we infer $X=0$. That means $\ker\, \sigma_{\textbf{S},\textbf{T}}=\{0\}$.
\end{proof}

\begin{thm}\label{main2}
Let $\text{\bf{T}}, \text{\bf{S}}\in \mathcal{B}^{m}_1(\Omega)$,
and $\text{\bf{T}}\sim_{u} (\text{\bf{M}}^*_z, \mathcal{H}_{K_{\text{\bf{T}}}}),\text{\bf{S}}\sim_{u} (\text{\bf{M}}^*_z, \mathcal{H}_{K_{\text{\bf{S}}}}).$
Then $\text{\bf{T}}\sim_{s}\text{\bf{S}}$ if and only if there exist two operator tuples $\widetilde{\text{\bf{T}}}=(T_{1}, \cdots, T_{m}), \widetilde{\text{\bf{S}}}=(S_{1}, \cdots, S_{m})\in \mathcal{B}^{m}_{2}(\Omega)$ such that
\begin{enumerate}
  \item [(1)]
 $T_{i}=\left (\begin{smallmatrix}
 T_{00}^i & T_{01}^i\\
 0 & T_{11}^i\\
\end{smallmatrix}\right ),S_{i}=\left(\begin{smallmatrix}
 S_{00}^i & S_{01}^i \\
 0 & S_{11}^i \\
\end{smallmatrix}\right),1\leq i\leq m$, where
$\text{\bf{T}}_{01}\in \mbox{ran}\,\sigma_{\text{\bf{T}}_{00},\text{\bf{T}}},\text{\bf{S}}_{01}\in \mbox{ran}\,\sigma_{\text{\bf{S}},\text{\bf{S}}_{11}}$ and $\mbox{ker}\ \sigma_{\text{\bf{T}}_{00},\text{\bf{S}}}=\mbox{ker}\ \sigma_{\text{\bf{S}}_{11},\text{\bf{T}}}=\{0\}$;
  \item [(2)] the bundles $E_{\widetilde{\text{\bf{T}}}}$ and $E_{\widetilde{\text{\bf{S}}}}$ of $\widetilde{\text{\bf{T}}}$ and $\widetilde{\text{\bf{S}}}$ are locally equivalent on an open dense subset of $\Omega$, the complement of the coalescing set for the curvature of $E_{\widetilde{\text{\bf{T}}}}$.  
\end{enumerate}
\end{thm}

\begin{proof}

From the proof of Theorem \ref{main9.26} and Corollary \ref{2110.26}, we only need to prove that there exists an $m$-tuple
$\textbf{S}_{11}\in \mathcal{B}^{m}_1(\Omega)$ satisfying the condition
$\mbox{ker}\ \sigma_{\textbf{S}_{11},\textbf{T}}=\{0\}$.

We first choose a generalized Bergman kernel $K_{\widehat{\textbf{S}}}$ on $\Omega\times\Omega$ (this concept was introduced by R.E. Curto and N. Salinas in \cite{CS2}),
which satisfies $\lim\limits_{dist(w,\partial \Omega)\rightarrow 0}K_{\widehat{\textbf{S}}}(w,w)=\infty$.
Set $K_{\textbf{S}_{11}}:=K_{\widehat{\textbf{S}}}\cdot K_{\textbf{T}}$. By \cite{CS2} and Theorem 2.6 in \cite{Salinas}, we know that $K_{\textbf{S}_{11}}(w,w)$ is also a generalized Bergman kernel and there exists
$\textbf{S}_{11}\in \mathcal{B}^{m}_1(\Omega)$ such that $\textbf{S}_{11}\sim_{u} (\textbf{M}^*_z, \mathcal{H}_{K_{\textbf{S}_{11}}})$.
Furthermore, we have that
$$\lim\limits_{dist(w,\partial \Omega)\rightarrow 0}\frac{K_{\textbf{T}}(w,w)}{K_{\textbf{S}_{11}}(w,w)}=\lim\limits_{dist(w,\partial \Omega)\rightarrow 0}\frac{1}{K_{\widehat{\textbf{S}}}(w,w)}=0.$$
By Lemma \ref{intertwining}, we know $\mbox{ker}\ \sigma_{\textbf{S}_{11},\textbf{T}}=\bigcap_{i=1}^{m}\mbox{ker}\ \sigma_{S_{11}^i,T_{11}^i}=\{0\}$.
This completes the proof.
\end{proof}

\subsection{Application}
Let $T$ be a bounded operator on some Hilbert $\mathcal{M}$, and $\mathcal{M}$ be a subspace of Hilbert space $\mathcal{N}$. A bounded operator $S$ on $\mathcal{N}$ is a dilation of $T$ if $P_{\mathcal{M}}S|_{\mathcal{M}}=T. $
We know that the adjoint of multiplication operator $M_z$ on Hardy space $\mathcal{H}$ is the Cowen-Douglas operator with index one over $\mathbb{D}$. Due to the fact $M_z\sim_s M_z|_{\mathcal K}$ for any invariant subspace $\mathcal{K}$ of $M_z$, we have that there exist plenty of operators such that their dilation is $M^*_z$ and they are all similar to $M^*_z$. Thus the following question is natural:

{\bf Question:}\,\,For any Cowen-Douglas operator $S$, is there a Cowen-Douglas operator $T$ such that $S$ is a dilation of $T$ but not similar to $T$?

By using the main theorem of this paper, we give lots of positive examples for this question.

In \cite{FJJ}, J.S. Fang, C.L. Jiang and the second author introduced an operator in the form of $T_{x}=\left(\begin{smallmatrix}
 T & x\otimes e_0 \\
 0 & M_{z}^{*} \\
\end{smallmatrix}\right),$ where $M_{z}^{*}$ is the adjoint of multiplication operator on Hardy space $\mathcal{H}$, $M_{z}^{*}e_0=0$
and $T\in \mathcal{L}(\mathcal{H})$ with spectral radius $r(T)<1$, $x\in \mathcal{H}$.
They proved that $T_{x}$ is a Cowen-Douglas operator with index one over $\Sigma$, a connected component of $\mathbb{D}\setminus\sigma(T)$ which contains $\{w\in \mathbb{D}: r(T)<|w|<1\}$.
Here  $T_{x} \in \mathcal{L}(\mathcal{H}\oplus \mathcal{H})$ is a dilation of $M_{z}^{*}$, since $P_{0\oplus \mathcal{H}}T_x|_{0\oplus \mathcal{H}}=M_{z}^{*}$.
In the following lemma, we replace the adjoint of the multiplication operator on Hardy space with a general Cowen-Douglas operator, and find the result is still valid.
\begin{lem}\label{2110.51}\cite{FJJ}
Let $S\in \mathcal{B}^{1}_{n}(\mathbb{D})\cap \mathcal{L}(\mathcal{H})$ and $S\sim_{u}(M_{z}^{*},\mathcal{H}_{K})$.
Suppose that $T\in \mathcal{L}(\mathcal{H})$ with spectral radius $r(T)<1$ and $T_{S,x}=\left(\begin{smallmatrix}
 T & x\otimes e_0 \\
 0 & S \\
\end{smallmatrix}\right),$
where $e_0=K(\cdot,0)\xi_{0}$ for some $\xi_0\in\mathbb{C}^{n}$, $x\in\mathcal{H}$.
Let $\Sigma$ be the connected component of $\mathbb{D}\setminus\sigma(T)$ which contains $\{w\in \mathbb{D}: r(T)<|w|<1\}$.
Then we have the following:
\begin{enumerate}
  \item For $w\in \mathbb{D}\setminus\sigma(T)$, $dim\ ker(T_{S,x}-w)=n$ and $ran(T_{S,x}-w)=\mathcal{H}\oplus\mathcal{H}$;
  \item $T_{S,x}\in \mathcal{B}^{1}_{n}(\Sigma)$ if and only if $x$ is a cyclic vector of $T$, i.e., $\bigvee_{n\geq0}\{T^{n}x\}=\mathcal{H}$.
\end{enumerate}
\end{lem}
\begin{proof}

In the following, we will describe $\ker(T_{S,x}-w)$ for $w\in \mathbb{D}\setminus\sigma(T)$. Suppose that $\left(\begin{smallmatrix}
 y_1  \\
 y_2 \\
\end{smallmatrix}\right)\in\ker(T_{S,x}-w), w\in \mathbb{D}\setminus\sigma(T)$. This is equivalent to $(T-w)y_1+(x\otimes e_0)y_2=0$ and $(S-w)y_2=0, w\in \mathbb{D}\setminus\sigma(T)$.
Without loss of generality, we assume that $S=(M_{z}^{*},\mathcal{H}_{K})$. Then $\ker(S-w)=\{K(\cdot,\bar{w})\xi,\xi\in\mathbb{C}^{n}\}$.
That means
$y_2=K(\cdot,\bar{w})\xi$ for some $\xi\in\mathbb{C}^{n}$. Note that $T-w$ is invertible when $w\in \mathbb{D}\setminus\sigma(T)$. We have
$$y_1=-(T-w)^{-1}(x\otimes e_0)y_2=-(T-w)^{-1}\langle K(0,\bar{w})\xi,\xi_0\rangle x. $$Thus, $-(T-w)^{-1}\langle K(0,\bar{w})\xi,\xi_0\rangle x\oplus K(\cdot,\bar{w})\xi$ is an eigenvector of $T_{S,x}$ with eigenvalue $w\in \mathbb{D}\setminus\sigma(T)$ for $\xi\in\mathbb{C}^{n}$. Since $dim \ker(S-w)=n$, we infer $dim\ ker(T_{S,x}-w)=n$.

For any $w\in \mathbb{D}\setminus\sigma(T)$, $T_{S,x}-w$ is surjective if for every $\left(\begin{smallmatrix}
 y_1'  \\
 y_2' \\
\end{smallmatrix}\right)\in \mathcal{H}\oplus\mathcal{H}$, there exists $\left(\begin{smallmatrix}
 y_1  \\
 y_2 \\
\end{smallmatrix}\right)\in \mathcal{H}\oplus\mathcal{H}$ such that $\left(\begin{smallmatrix}
 T-w & x\otimes e_0 \\
 0 & S-w \\
\end{smallmatrix}\right)\left(\begin{smallmatrix}
 y_1  \\
 y_2 \\
\end{smallmatrix}\right)=\left(\begin{smallmatrix}
 (T-w)y_1+(x\otimes e_0)y_2 \\
  (S-w)y_2 \\
\end{smallmatrix}\right)=\left(\begin{smallmatrix}
 y_1'  \\
 y_2' \\
\end{smallmatrix}\right)$.
The existence of $y_2$ is clear, since $S$ is a Cowen-Douglas operator. If we take $y_1=(T-w)^{-1}(y_1'-(x\otimes e_0)y_2)$, then the last equation holds. This proves the statement (1).

In order to get statement (2), by (1), we only need to prove that $$\bigvee_{w\in\Sigma}\{-(T-w)^{-1}\langle K(0,\bar{w})\xi,\xi_0\rangle x\oplus K(\cdot,\bar{w})\xi,\xi\in\mathbb{C}^{n}\}=\mathcal{H}\oplus\mathcal{H}$$ and $\bigvee_{n\geq0}\{T^{n}x\}=\mathcal{H}$ are equivalent.
Suppose that there exists an $x_1\oplus x_2\in \mathcal{H}\oplus\mathcal{H}$ such that $\big\langle -(T-w)^{-1}\langle K(0,\bar{w})\xi,\xi_0\rangle x\oplus K(\cdot,\bar{w})\xi, x_1\oplus x_2\big\rangle=0$. Then
$$\big\langle (w-T)^{-1}\langle K(0,\bar{w})\xi,\xi_0\rangle x,x_1\big\rangle=-\langle K(\cdot,\bar{w})\xi, x_2\rangle.$$
Note that $K(\cdot,\bar{w})\xi$ is analytic on $\mathbb{D}$, $(w-T)^{-1}=\frac{1}{w}\sum_{n=0}^\infty(\frac{T}{w})^n$ for $|w|>r(T)$ and
$\langle (w-T)^{-1}\langle K(0,\bar{w})\xi,\xi_0\rangle x,x_1\big\rangle$ is analytic when $|w|>r(T)$.
Then
$\big\langle (w-T)^{-1}\langle K(0,\bar{w})\xi,\xi_0\rangle x,x_1\big\rangle=-\langle K(\cdot,\bar{w})\xi, x_2\rangle,r(T)<|w|<1$.
Thus, by the analytic continuation theorem, we know that $\big\langle (w-T)^{-1}\langle K(0,\bar{w})\xi,\xi_0\rangle x,x_1\big\rangle$ is analytic on $\mathbb{C}$. Since $\lim_{|w|\rightarrow\infty}\langle \sum_{n=0}^\infty(\frac{T^nx}{w^{n+1}}),x_1\rangle=0$,
$\big\langle (w-T)^{-1}\langle K(0,\bar{w})\xi,\xi_0\rangle x,x_1\big\rangle$ is a bounded entire function on $\mathbb{C}$. Then
$\sum_{n=0}^\infty\langle T^nx,x_1\rangle \frac{1}{w^{n+1}}=0.$ Therefore, $\langle T^nx,x_1\rangle=0,n\geq0$.
Suppose $x$ is a cyclic vector of $T$. This implies
$x_1=0$. From $\bigvee_{w\in\Sigma}\{ K(\cdot,\bar{w})\xi,\xi\in\mathbb{C}^{n}\}=\mathcal{H}$, we see that $x_2=0$.
That means $\bigvee_{w\in\Sigma}\ker(T_{S,x}-w)=\mathcal{H}\oplus\mathcal{H}.$  Suppose $x$ is not a cyclic vector of $T$.
Let $0\neq x_1\perp \{T^nx: n\geq0\}$. Then $(x_1\oplus0)\perp \bigvee_{w\in\Sigma}\ker(T_{S,x}-w)$ and therefore
$\bigvee_{w\in\Sigma}\ker(T_{S,x}-w)\neq\mathcal{H}\oplus\mathcal{H}$. This completes the proof.
\end{proof}

\begin{prop}
Let $T_{S,x}$ be the operator in  Lemma \ref{2110.51} and $T\in \mathcal{L}(\mathcal{H})$ with spectral radius $r(T)<1$.
Suppose that $\lim_{|w|\rightarrow r(T)}\|(T-w)^{-1}x\|=\infty$ and $x$ is a cyclic vector of $T$, then $T_{S,x}$ is not similar to $S$.
\end{prop}
\begin{proof}
Suppose that $T_{S,x}$ is similar to $S$. By Theorem \ref{main9.26}, without losing generality, there exists a bounded linear operator $X$ such that $Y=(I+X^*X)^{\frac{1}{2}}$ and $T_{S,x}Y=YS$.
By Lemma \ref{2110.51} and $x$ is a cyclic vector of $T$, we have that $T_{S,x}\in \mathcal{B}^{1}_{n}(\Sigma)$, where $\Sigma$ is the connected component of $\mathbb{D}\setminus\sigma(T)$ which contains $\{w\in \mathbb{D}: r(T)<|w|<1\}$.
Note that $-(T-w)^{-1}\langle K(\cdot,\bar{w})\xi,e_0\rangle x\oplus K(\cdot,\bar{w})\xi\in \ker(T_{S,x}-w)$ for $\xi\in \mathbb{C}^n$ and $w\in\Sigma$.
By Proposition 2.4 in \cite{R} and $T_{S,x}Y=YS$, we find $\xi_w\in\mathbb{C}^n$ such that
$$(I+X^*X)^{\frac{1}{2}}K(\cdot,\bar{w})\xi_{w}=-(T-w)^{-1}\langle K(0,\bar{w})\xi_w,\xi_0\rangle x\oplus K(\cdot,\bar{w})\xi_{w}.$$
This implies that $\|(I+X^*X)^{\frac{1}{2}}K(\cdot,\bar{w})\xi_{w}\|^{2}\leq(1+\|X\|^{2})\langle K(\bar{w},\bar{w})\xi_w,\xi_w\rangle$ and
$$\begin{array}{lll}
\|(I+X^*X)^{\frac{1}{2}}K(\cdot,\bar{w})\xi_{w}\|^{2}
&=&\|(T-w)^{-1}\langle K(0,\bar{w})\xi_w,\xi_0\rangle x\|^{2}+\langle K(\bar{w},\bar{w})\xi_w,\xi_w\rangle\\
&=&|\langle K(0,\bar{w})\xi_w,\xi_0\rangle |^{2}\|(T-w)^{-1}x\|^{2}+\langle K(\bar{w},\bar{w})\xi_w,\xi_w\rangle.
\end{array}$$
Then we have  
$0\leq|\langle K(0,\bar{w})\xi_w,\xi_0\rangle|^2\|(T-w)^{-1}x\|^{2}\leq \|X\|^{2}\langle K(\bar{w},\bar{w})\xi_w,\xi_w\rangle$.
Note that $\langle K(\bar{w},\bar{w})\xi_w,\xi_w\rangle,w\in\mathbb{D}$ and $X$ are bounded, then
$|\langle K(0,\bar{w})\xi_w,\xi_0\rangle |^{2}\|(T-w)^{-1}x\|^{2}$ is bounded.
Since $\lim_{|w|\rightarrow r(T)}\|(T-w)^{-1}x\|=\infty, $
$|\langle K(0,\bar{w})\xi_w,\xi_0\rangle |\to 0$ when $|w|\rightarrow r(T)$. We know that $\langle K(0,\bar{w})\xi_w,\xi_0\rangle$ is holomorphic,
by the maximum modulus principle of holomorphic function, then $\langle K(0,\bar{w})\xi_w,\xi_0\rangle=0$ for all $w\in\Sigma$.
This is a contradiction. So $T_{S,x}$ is not similar to $S$.
\end{proof}


Upon using our main theorem (Theorem \ref{main9.26}), a new proof of the sufficiency of A.L. Shields' similarity theorem in \cite{SAL} is given.
\begin{ex}\label{ALex}
Let $\mathbf{T}=(T_{11}^1, \cdots, T_{11}^m),\mathbf{S}=(S_{00}^1, \cdots, S_{00}^m)\in \mathcal{B}_{1}^{m}(\Omega)$ and $\mathbf{T}\sim_u(\mathbf{M}_z^*,\mathcal{H}_{K_0})$, $\mathbf{S}\sim_u(\mathbf{M}_z^*,\mathcal{H}_{K_1})$, where $K_i(z,w)=\sum_{\alpha\in\mathbf{Z}_+^m} \rho_i(\alpha)z^\alpha\bar{w}^\alpha$ and $\rho_i(\alpha)>0$ for $i=0,1,\alpha\in\mathbf{Z}_+^m$.
If $\frac{\rho_1(\alpha)}{\rho_0(\alpha)}$ is bounded for all $\alpha\in\mathbf{Z}_+^m$, then $\mathbf{T}\sim_s \mathbf{S}$.
\end{ex}
\begin{proof}
By section 2.1 of \cite{LQ} due to Q. Lin and Theorem \ref{main2} in this paper, we find that there exist tuples $\mathbf{T}_0,\mathbf{S}_1$ such that $\ker\,\sigma_{\mathbf{T}_0,\mathbf{S}}=\ker\,\sigma_{\mathbf{S}_1,\mathbf{T}}=\{0\}$.
Since $\frac{\rho_1(\alpha)}{\rho_0(\alpha)}$ is bounded for all $\alpha\in\mathbf{Z}_+^m$, it follows that there exist $l>0$ such that $l^2\rho_1(\alpha)-\rho_0(\alpha)>0$ and $b_\alpha:=\sqrt{l^2\frac{\rho_1(\alpha)}{\rho_0(\alpha)}-1}$ is bounded for all $\alpha\in\mathbf{Z}_+^m$.
Select the non-zero holomorphic sections of $\mathbf{T}$ and $\mathbf{S}$ as $K_0(\cdot,\bar{w})$ and $lK_1(\cdot,\bar{w})$, respectively.
Define operator
$X:\mathcal{H}_{K_0}\rightarrow\mathcal{H}_{K_1}$ as $Xw^\alpha=b_\alpha w^\alpha$ for all $\alpha\in\mathbf{Z}_+^m$.
Thus, $l^2K_1(\bar{w},\bar{w})=K_0(\bar{w},\bar{w})+\|XK_0(\cdot,\bar{w})\|^2=\|(1+X^*X)^{\frac{1}{2}}K_0(\cdot,\bar{w})\|^2$.
Note that $I+X^*X$ and $I+XX^*$ are positive and invertible, there exist invertible operators $U_{01},U_{10}$ such that
$(1+XX^*)^{-1}=U^*_{10}U_{10}, (1+X^*X)^{-1}=U^*_{01}U_{01}$. Choosing $\widetilde{Y}\in\ker\sigma_{\mathbf{S},\mathbf{S}_1}$. Set $Y=\widetilde{Y}+U_{01}X^*U_{10}^{-1}$,
$\widetilde{\textbf{T}}=(T_{1}, \cdots, T_{m}), \widetilde{\textbf{S}}=(S_{1}, \cdots, S_{m})$ and $U=\left(\begin{smallmatrix}
 U_{01}X^* & U_{01}\\
 U_{10} & -U_{10}X \\
 \end{smallmatrix}\right)$, where $T_{i}=\left (\begin{smallmatrix}
 T_{00}^i & T_{01}^i\\
 0 & T_{11}^i\\
\end{smallmatrix}\right ),S_{i}=\left(\begin{smallmatrix}
 S_{00}^i & S_{01}^i \\
 0 & S_{11}^i \\
\end{smallmatrix}\right),1\leq i\leq m,$ and $\textbf{T}_{01}= -\sigma_{\textbf{T}_{00},\textbf{T}}(X),\textbf{S}_{01}=-\sigma_{\textbf{S},\textbf{S}_{11}}(Y)$. By Theorem \ref{main9.26}, we know that $U$ is  unitary and $U\widetilde{\textbf{T}}=\widetilde{\textbf{S}}U$. Thus we have $\mathbf{T}\sim_{s}\mathbf{S}$.
\end{proof}


\begin{ex}
Let $T_i\in \mathcal{B}_{1}^{1}(\mathbb{D})$, $T_i\sim_u(M_z^*,\mathcal{H}_{K_i}),i=0,1$, where $K_0(z,w)$ is the reproducing kernel in Example \ref{ALex}. 
If $K_1(z,w)-K_0(z,w)=P(z,\bar{w})$ is a polynomial and positive over $\mathbb{D}\times \mathbb{D}$, then $T_0\sim_s T_1$.
\end{ex}
\begin{proof}
Since the symmetry of $K_i,i=0,1$, $P$ is also symmetric. Without losing generality, we assume that $P(z,\bar{w})=\sum_{p,q=0}^ma_{pq}z^p\bar{w}^q,$ $z,\bar{w}\in\mathbb{D}$ for some positive $m$.
This implies that matrix $A:=(a_{pq})_{p,q=0}^{m}$ is positive. By diagonalization of $A$, there exist $\{\phi_i\}_{i=0}^m\subset H^{\infty}(\mathbb{D})$, a
set of bounded holomorphic function on $\mathbb{D}$, such that $P(z,\bar{w})=\sum_{j=0}^m\phi_j(z)\overline{\phi_j(w)}, z,\bar{w}\in \mathbb{D}$.
Let $K_0(z,w)=\sum_{n=0}^\infty e_n(z)e_n^*(w)$ for the orthogonal normal basis $\{e_n\}_{n=0}^\infty$ of $\mathcal{H}_{K_0}$.
Set $\phi_j(z)=\sum_{i=0}^mb_{ji}e_i(z)$, $b_{ji}\in \mathbb{C},0\leq j\leq m. $
A linear operator $X$ is defined as the
following:
$$X(e_{i}):=\left\{
\begin{array}{cc}
\sum_{j=0}^{m}b_{ji}e_j,\ 0\leq i\leq m;\\
0,\qquad\ i>m.
\end{array}
\right.
$$
Then
$$\begin{array}{lll}
\|X\|&=&\sup\limits_{\|y\|=1}\big\|Xy\big\|=\sup\limits_{\|y\|=1}\Big\|X\sum\limits_{n=0}^{\infty}b_ne_n\Big\|\\&=&\sup\limits_{\|y\|=1}\Big\|\sum\limits_{n=0}^{m}b_n\Big(\sum\limits_{j=0}^{m}b_{jn}e_j\Big)\Big\|\\
&=&\sup\limits_{\|y\|=1}\Big\|\sum\limits_{j=0}^{m}\sum\limits_{n=0}^{m}b_nb_{jn}e_j\Big\|\\
&\leq&\Big(\sum\limits_{n=0}^{m}b^{2}_n\Big)^{\frac{1}{2}}\sum\limits_{j=0}^{m}\Big(\Big\|\sum\limits_{n=0}^{m}b_{jn}^{2}\Big\|^{\frac{1}{2}}\Big)\leq (m+1)M,
\end{array}$$
where $y=\sum_{n=0}^{\infty}b_ne_n\in\mathcal{H}$, $M=\max\limits_{0\leq i\leq m}\{\|\phi_i\|\}$. Thus, $X$ is bounded. Note that
$XK_0(z,\bar{w})
=\sum\limits_{i=0}^{\infty}Xe_i(w)e_i(z)
=\sum\limits_{i=0}^{m}\sum\limits_{j=0}^{m}b_{ji}e_j(w)e_i(z)=\sum\limits_{j=0}^{m}\sum\limits_{i=0}^{m}b_{ji}e_i(z)e_j(w)
=\sum\limits_{j=0}^{m}\phi_j(z)e_{j}(w),$
we then have $\|XK_0(z,\bar{w})\|^2=\sum\limits_{i=0}^{m}|\phi_i(z)|^2$ and $\|XK_0(\cdot,\bar{w})\|^2=P(w,\bar{w})$.
Similarly to the proof of Example \ref{ALex}, we deduce that $T_0\sim_{s}T_1$.
\end{proof}

\begin{ex}\label{ex3.12}
Let $\textbf{T}_i=(T_{i_1},\cdots,T_{i_m})\in\mathcal{B}_1^m(\Omega)\cap \mathcal{L}(\mathcal{H}_i)$ and $\textbf{T}_i\sim_{u}(\textbf{M}_z^*,\mathcal{H}_{K_i}),i=0,1$.
If there exists a uniformly bounded positive sequence
$\{\lambda_\alpha\}_{\alpha\in\mathbf{Z}_+^m}$ such that $K_1(z,w)-K_0(z,w)=\sum\limits_{\alpha\in\mathbf{Z}_+^m}\lambda_\alpha e_\alpha(z)e^*_\alpha(w)$
for some orthonormal basis $\{e_{\alpha}\}_{\alpha\in\mathbf{Z}_+^m}$ of $\mathcal{H}_0$,
then $\textbf{T}_0\sim_{s}\textbf{T}_1$.
\end{ex}
\begin{proof}
Without loss of generality, we assume that $\textbf{T}_i=(\textbf{M}_z^*,\mathcal{H}_{K_i}),i=0,1$.
Since $\textbf{T}_i\in\mathcal{B}_1^m(\Omega)$, we know that $K_i(\cdot,\bar{w})\in\ker(\textbf{T}_i-w)$ for $w\in\Omega$ and $i=0,1$. Then
\begin{equation}\label{2111.51}
\begin{array}{lll}
K_1(w,w)-K_0(w,w)
&=&\sum\limits_{\alpha\in\mathbf{Z}_+^m}\lambda_\alpha|e_\alpha(w)|^2\\
&=&\sum\limits_{\alpha\in\mathbf{Z}_+^m}\lambda_\alpha|\langle K_0(\cdot,\bar{w}),e_\alpha\rangle|^2\\
&=&\sum\limits_{\alpha\in\mathbf{Z}_+^m}\lambda_\alpha\langle K_0(\cdot,\bar{w}),e_\alpha\rangle \langle e_\alpha,K_0(\cdot,\bar{w})\rangle\\
&=& \Big\langle \sum\limits_{\alpha\in\mathbf{Z}_+^m}\lambda_\alpha\langle K_0(\cdot,\bar{w}),e_\alpha\rangle e_\alpha,K_0(\cdot,\bar{w})\Big\rangle.
\end{array}
\end{equation}
Let $X_1:=\sum_{\alpha\in\mathbf{Z}_+^m}\lambda_\alpha e_\alpha\otimes e_\alpha$. 
Further, equation (\ref{2111.51}) can be written as $K_1(w,w)-K_0(w,w)=\langle X_1K_0(\cdot,\bar{w}),K_0(\cdot,\bar{w})\rangle. $
It is easy to see that $X_1$ is positive, since $\lambda_\alpha>0$ for all $\alpha\in\mathbf{Z}_+^m$. Thus, there exists an operator $X_2$ such that $X_1=X_2^*X_2$ and $K_1(w,w)-K_0(w,w)=\|X_2K_0(\cdot,\bar{w})\|^2.$
Similarly to the proof of Example \ref{ALex}, we deduce that $\mathbf{T}_0\sim_{s}\mathbf{T}_1$.
\end{proof}

\begin{ex}
Let $\textbf{T}_i=(T_{i_1},\cdots,T_{i_m})\in\mathcal{B}_1^m(\Omega)\cap \mathcal{L}(\mathcal{H}_i)$ and $\textbf{T}_i\sim_{u}(\textbf{M}_z^*,\mathcal{H}_{K_i}),i=0,1$.
Suppose that $L^2(X,\mu)$ is separable for some $\sigma$-finite measure space $(X, \mu)$.
If there exists $\phi\in L^\infty(X,\mu)$ (need not to be holomorphic) such that $K_1(w,w)=(1+|\phi(w)|^2)K_0(w,w)$, then $\textbf{T}_0\sim_{s}\textbf{T}_1$.
\end{ex}
\begin{proof}
From $\phi\in L^\infty(X,\mu)$ for some $\sigma$-finite measure space $(X, \mu)$, then there is a multiplication operator $M_\phi: L^2(X,\mu)\rightarrow L^2(X,\mu)$ is defined by $M_\phi f(x)=\phi(x)f(x)$, and it satisfies
$\|M_\phi f\|_2=(\int_X|\phi f|^2d\mu)^\frac{1}{2}\leq(\int_X(\|\phi\|_\infty |f|)^2d\mu)^\frac{1}{2}\leq \|\phi\|_\infty \|f\|,f\in L^2(X,\mu)$,
thus $M_\phi$ is bounded. 
For $f,g\in L^2(X,\mu)$, we have
$\langle M_\phi f, g\rangle=\int_X (\phi f)\bar{g}d\mu=\int_X f\overline{(\bar{\phi} g)}d\mu=\langle f, M_{\bar{\phi}}g\rangle,$
which implies $M_\phi^*=M_{\bar{\phi}}$.

From $K_1(w,w)=(1+|\phi(w)|^2)K_0(w,w)$, we have
$\|K_1(\cdot,\bar{w})\|^2=(1+|\phi(w)|^2)\|K_0(\cdot,\bar{w})\|^2$. 
Since $L^2(X,\mu)$ is separable, there is a unitary operator $U: \mathcal{H}_{K_0}\rightarrow L^2(X,\mu)$ such that
$$\begin{array}{lll}
\|K_1(\cdot,\bar{w})\|^2
&=&(1+|\phi(w)|^2)\|UK_0(\cdot,\bar{w})\|^2\\
&=&\langle (1+|\phi(w)|^2)UK_0(\cdot,\bar{w}),UK_0(\cdot,\bar{w})\rangle\\
&=&\langle (I+M_\phi^*M_\phi)UK_0(\cdot,\bar{w}),UK_0(\cdot,\bar{w})\rangle\\
&=&\|(I+U^*M_\phi^*M_\phi U)^\frac{1}{2}K_0(\cdot,\bar{w})\|^2.
\end{array}$$
Without loss of generality, we assume that $\textbf{T}_i=(\textbf{M}_z^*,\mathcal{H}_{K_i})$, then $K_i(\cdot,\bar{w})\in\ker(\textbf{T}_i-w)$ for $w\in\Omega$ and $i=0,1$.
Similarly to the proof of Example \ref{ALex}, we deduce that $\mathbf{T}_0\sim_{s}\mathbf{T}_1$.
\end{proof}

\section{A Subclass $N\mathcal{FB}_{n_0,n_1}^{m}(\Omega)$ of Cowen-Douglas tuples }

Let $M\ddot{o}b$ denote the group of all biholomorphic automorphisms of $\mathbb{D}$.
Recall that a bounded operator $T$ is said to be homogeneous if the spectrum
$\sigma(T)$ of $T$ is contained in $\overline{\mathbb{D}}$ and for every $\phi\in M\ddot{o}b$, $\phi(T)$ is unitarily equivalent to $T$.
The concept of homogeneous operator can be extended to the commuting operator tuple.
When $\mathcal{D}$ is a bounded symmetric domain, an $m$-tuple
$\textbf{T}=(T_1,\cdots,T_m)$ of commuting bounded operators is said to be homogeneous with respect
to $G$ if
their joint Taylor spectrum is contained in $\overline{\mathcal{D}}$ and for every holomorphic automorphism $\phi\in G$, $\phi(\textbf{T})$ is unitarily equivalent to $\textbf{T}$ (see \cite{BM3,MS3}).
The topic of homogeneous operators and tuples received much attention \cite{BM1,BM2,BM3,Misra,Hazra,AM,KM0,SS}, mostly by using representation theory of Lie groups and complex geometry.
G. Misra in \cite{Misra} has fully characterized the homogeneous operators in $\mathcal{B}_1^1(\mathbb{D})$, e.g. proving that for any
$j>0$, the unilateral shift operator with weight sequence $\{\sqrt{\frac{i+1}{i+j}}\}_{i=0}^\infty$ is homogeneous.
The homogeneous operator in $\mathcal{B}_1^1(\mathbb{D})$ not only provides us with a model, but also helps us to study the properties and similarity of other operators.
In \cite{AM}, A. Koranyi and G. Misra completed the classification of irreducible homogeneous operators in $\mathcal{B}_n^1(\mathbb{D})$.

For $\alpha\in\mathbb{C}$ and homogeneous operators $T_0,T_1$ acting on $\mathcal{H}$, $\left(\begin{smallmatrix}
 T_{0} & \alpha(T_{0}-T_{1}) \\
 0 & T_{1} \\
 \end{smallmatrix}\right)$ is homogeneous if $T_0$ and $T_1$ have the same associated unitary representation given by A. Koranyi in Lemma 2.1 \cite{AK}.
Let $t_n=t_n(a,b)=\sqrt{\frac{n+a}{n+b}},n\in\mathbb{Z}$ for $a,b\in(0,1),a\neq b$.
Define the operator $T=T(a,b)$ as $Te_n=t_ne_{n+1}$ for the natural basis $\{e_n\}_{n\in\mathbb{Z}}$ of $l^2(\mathbb{Z})$. It is shown to be homogeneous.
Let $\alpha>0$. Defining $\tilde{T}=\tilde{T}(a,b,\alpha)=\left(\begin{smallmatrix}
 T(a,b) & \alpha(T(a,b)-T(b,a)) \\
 0 & T(b,a) \\
 \end{smallmatrix}\right)$ on $\mathcal{H}\oplus\mathcal{H}$ and rearranging the bases, we obtain a block matrix so that it is
$\tilde{T_n}=\tilde{T_n}(a,b,\alpha)$ at $(n+1,n)$-position and the rest are 0.
It is shown that $\tilde{T}$ is homogeneous and irreducible in \cite{AK}.
This is the first example of irreducible bi-lateral homogeneous 2-shifts with three parameters due to A. Koranyi.
Next, another bi-lateral homogeneous 2-shift introduced by S. Hazra in \cite{Hazra}.
Let $B(s)$ and $B$ be bi-lateral shifts, the weight sequence of $B(s)$ be $w_n=\frac{n+\frac{1+\lambda}{2}+s}{n+\frac{1+\lambda}{2}-s}(s\neq0)$, and $B$ be unweighted.
Then operators $B(s)$ and $B$ are homogeneous in Theorem 5.2 of \cite{BM2}. For $\alpha>0$, define $B(\lambda,s,\alpha)=\left(\begin{smallmatrix}
 B(s) & \alpha(B(s)-B) \\
 0 & B \\
 \end{smallmatrix}\right)$. It is also homogeneous by Lemma 2.1 of \cite{AK}.
It is proved in \cite{Hazra} that $B(\lambda,s,\alpha)$ is irreducible and the homogeneous operators defined by A. Koranyi and S. Hazra, respectively, are mutually unitarily inequivalent.

Inspired by the above results, we here define a new class of tuples of commuting bounded operators. With the help of this class, we discuss the similarity of commuting tuples.
In addition, we know that the Cowen-Douglas class is a very rich operator and tuple class, including many homogeneous operators, normal operators and so on. The structure of the elements in $\mathcal{B}_{n}^m(\Omega)$ is very complicated, so that we still cannot clearly describe their similarity. Therefore, it is necessary to investigate a subclass of $\mathcal{B}_{n}^m(\Omega)$.


\subsection{Definitions}
In what follows, we assume that $n_0,n_1$ are positive integers.
\begin{defn}\label{maindefn}
Let $\textbf{T}_{ii}=(T_{ii}^1, \cdots, T_{ii}^m)\in \mathcal{B}_{n_i}^{m}(\Omega),i=0,1$ and
$\textbf{T}_{01}=(T_{01}^1, \cdots, T_{01}^m)$ be a commuting $m$-tuple of bounded operators.
Suppose that the $m$-tuple $\textbf{T}=(T_{1}, \cdots, T_{m})$ satisfies $T_{j}=\left(\begin{smallmatrix}
 T_{00}^j & T_{01}^j \\
 0 & T_{11}^j \\
 \end{smallmatrix}\right)$ for $1\leq j\leq m$.
We call $\textbf{T}\in N\mathcal{FB}_{n_0,n_1}^{m}(\Omega)$, if $\textbf{T}_{01}\in \mbox{ran}\,\sigma_{\textbf{T}_{00},\textbf{T}_{11}}$.
\end{defn}
By Corollary \ref{2110.26}, we see that tuples in $N\mathcal{FB}_{n_0,n_1}^{m}(\Omega)$ are Cowen-Douglas tuples with index $n_0+n_1$ over $\Omega$.
If $n_0=n_1=n$, the class $N\mathcal{FB}_{n_0,n_1}^{m}(\Omega)$ can be expressed as $N\mathcal{FB}_{2n}^{m}(\Omega)$.
\begin{rem}
Suppose that tuple $\textbf{T}$ satisfies the conditions of Definition \ref{maindefn} and there exists an operator $X$ such that $\textbf{T}_{01}=\sigma_{\textbf{T}_{00},\textbf{T}_{11}}(-X)$.
Then $\textbf{T}_{ii}$, $i=0,1$ are commuting $m$-tuples, which means $\textbf{T}$ is a commuting $m$-tuple, since
$$\begin{array}{lll}
T_{p}T_{q}&=&\begin{pmatrix}
 T_{00}^pT_{00}^q & -T_{00}^pT_{00}^qX+XT_{11}^pT_{11}^q \\
 0 & T_{11}^pT_{11}^q \\
 \end{pmatrix}\\
&=&\begin{pmatrix}
  T_{00}^qT_{00}^p & -T_{00}^qT_{00}^pX+XT_{11}^qT_{11}^p \\
 0 & T_{11}^qT_{11}^p \\
 \end{pmatrix}\\
&=&T_{q}T_{p}\,\ \mbox{\ for\, all}\  1\leq p, q\leq m.
\end{array}$$
\end{rem}
Based on Theorem 1.49 in \cite{JW}, C. Jiang, D.K. Keshari, G. Misra and the second author introduced an operator class, denoted by  $\mathcal{FB}^1_{n}(\Omega)$ in \cite{JJKMCR,JJKM}, which is norm dense in $\mathcal{B}^{1}_{n}(\Omega)$.
They also showed that the complete unitary invariants of operators in $\mathcal{FB}^1_{n}(\Omega)$ include the curvatures and the second fundamental forms of the diagonal operators. 
We will give the commuting tuple version of this kind of operator.
\begin{defn}\label{220418.1}
Let $\textbf{T}_{ii}=(T_{ii}^1, \cdots, T_{ii}^m)\in \mathcal{B}_{1}^{m}(\Omega)\cap\mathcal{L}(\mathcal{H}_{i})^{m},i=0,1$.
Suppose that there exists $T_{01}\in \mathcal{L}(\mathcal{H}_{1},\mathcal{H}_{0})$ such that
$\textbf{T}=(T_{1}, \cdots, T_{m})$ is a commuting $m$-tuple with $T_{j}=\left(\begin{smallmatrix}
 T_{00}^i & T_{01}\\
 0 & T_{11}^i \\
 \end{smallmatrix}\right),\,1\leq j\leq m$.
We call $\textbf{T}\in\mathcal{FB}_{2}^{m}(\Omega)$, if $T_{01}\in \mbox{ker}\,\sigma_{\textbf{T}_{00},\textbf{T}_{11}}$.
\end{defn}
In order to show that the tuples in $N\mathcal{FB}^{m}_2(\Omega)$ may not belong to $\mathcal{FB}^{m}_2(\Omega)$, we need to introduce the following concept which is first
defined in \cite{JJ}.

\begin{defn}\cite{JJ}\label{ph} {\bf{Property (H)}}
Let $T_{ii}\in \mathcal{L}(\mathcal{H}_{i}),i=0,1$ and $T=\left(\begin{smallmatrix}
 T_{00} & XT_{11}-T_{00}X\\
 0 & T_{11} \\
 \end{smallmatrix}\right)\in \mathcal{B}^{1}_2(\Omega)$. We call $T$ satisfies the Property (H) if and only if the following statements hold:
If $Y\in \mathcal{L}(\mathcal{H}_{1}, \mathcal{H}_0)$ satisfies
\begin{enumerate}
\item[(i)] $T_{00}Y=YT_{11},$
\item[(ii)]$Y=T_{00}Z-ZT_{11}$ for some $Z$.
\end{enumerate}
Then $Y=0$.  That is equivalent to $\mbox{ker}\ \sigma_{T_{00},T_{11}}\cap \mbox{ran}\ \sigma_{T_{00},T_{11}}=\{0\}$.
\end{defn}

By Definition \ref{ph}, we see if  $T=\left(\begin{smallmatrix}
T_{00} & XT_{11}-T_{00}X\\
0 & T_{11} \\
\end{smallmatrix}\right)$ satisfies the Property $(H)$, and $XT_{11}\neq T_{00}X$, then $T$ does not belong to $\mathcal{FB}^{1}_2(\Omega)$. Otherwise,
$XT_{11}-T_{00}X\in \mbox{ker}\ \sigma_{T_{00},T_{11}}\cap \mbox{ran}\ \sigma_{T_{00},T_{11}}.$ That means  $XT_{11}=T_{00}X$. It is a contradiction.
In the following, we will give two results to show when $T$ would satisfy the Property $(H)$.
\begin{prop}\label{JiH}\cite{JJ} Let $T_0,T_1\in \mathcal{L}(\mathcal {H})$ and $S_1$ be the right inverse of $T_1$. If $\lim\limits_{n\rightarrow \infty} \frac{\|T^n_0\|\cdot \|S^n_1\|}{n}=0$, then the Property (H) holds, i.e. if there exists $X\in \mathcal{L}(\mathcal {H})$ such that
$T_0X=XT_1$ and $X=T_0Y-YT_1$, then $X=0$ (i.e. $\mbox{ker}\ \sigma_{T_0,T_1}\cap \mbox{ran}\ \sigma_{T_0,T_1}=\{0\}$).
\end{prop}
\begin{ex}\cite{JJ} Let $A,B\in \mathcal{B}^{1}_1(\mathbb{D})$ be backward shift operators with weighted sequences $\{a_i\}^{\infty}_{i=1}$ and $\{b_i\}^{\infty}_{i=1}$.
If $\lim\limits_{n\rightarrow \infty}n\frac{\prod\limits^n_{k=1}b_k}{\prod\limits^n_{k=1}a_k}=\infty$, then $ker\ \sigma_{A,B}\cap ran\ \sigma_{A,B}=\{0\}$.
\end{ex}
In \cite{JJKM}, it is proved that the unitary operator of intertwining two operators $T$ and $\tilde{T}$ in $\mathcal{FB}^{1}_n(\Omega)$ should be a
diagonal matrix.
From the proof of Lemma \ref{mainlemma}, it can be seen that the unitary operator intertwines the two tuples in the class $N\mathcal{FB}^{m}_{n}(\Omega)$ could be non-diagonal.
This is another reason why we study this new class.
Although the structures of tuples in the classes $N\mathcal{FB}_{2}^{m}(\Omega)$ and $\mathcal{FB}_{2}^{m}(\Omega)$ are quite different, the following proposition shows that they are also closely related. The unitary equivalence of the tuples in $ N\mathcal{FB}_{2}^{m}(\Omega)$ can always be related to the similarity of the tuples in $\mathcal{FB}_{2}^{m}(\Omega)$.

\begin{prop}\label{main1}
For $i=0,1,$ let $\text{\bf{T}}_{ii},\text{\bf{S}}_{ii}\in\mathcal{B}_{1}^{m}(\Omega), \text{\bf{T}}_{01},\text{\bf{S}}_{01}\in\mathcal{L}(\mathcal{H})^m$. Let
$\widetilde{\text{\bf{T}}}=(T_{1},\cdots, T_{m})$, $\widetilde{\text{\bf{S}}}=(S_{1}, \cdots, S_{m})\in N\mathcal{FB}_{2}^{m}(\Omega)$ with
$T_j=\left(\begin{smallmatrix}T_{00}^j&T_{01}^j\\0&T_{11}^j\end{smallmatrix}\right)$, $S_j=\left(\begin{smallmatrix}S_{00}^j&S_{01}^j\\0&S_{11}^j\end{smallmatrix}\right), 1\leq j\leq m$.
Suppose that $\mbox{ker}\, \sigma_{\text{\bf{T}}_{00},\text{\bf{S}}_{00}}=\mbox{ker}\, \sigma_{\text{\bf{S}}_{11},\text{\bf{T}}_{11}}=\{0\}$.
If $\widetilde{\text{\bf{T}}}\sim_{u}\widetilde{\text{\bf{S}}}$, then there exist operators $S_0,S_1$ and tuples $\widehat{\text{\bf{T}}}=(\widehat{T}_1,\cdots,\widehat{T}_m),\widehat{\text{\bf{S}}}=(\widehat{S}_1,\cdots,\widehat{S}_m)$ with
$\widehat{T}_i=\left(\begin{smallmatrix}
 S_{00}^i & S_0 \\
 0 & T_{00}^i\\
 \end{smallmatrix}\right),\widehat{S}_i=\left(\begin{smallmatrix}
 T_{11}^i & S_1 \\
 0 & S_{11}^i\\
 \end{smallmatrix}\right), 1\leq i\leq m$, such that $\widehat{\text{\bf{T}}},\widehat{\text{\bf{S}}}\in \mathcal{FB}^{m}_2(\Omega)$ and
$\widehat{\text{\bf{T}}}\sim_s \widehat{\text{\bf{S}}}$.
\end{prop}

\begin{proof}
Suppose that there exist $X,Y$ such that
$\textbf{T}_{01}=\sigma_{\textbf{T}_{00},\textbf{T}_{11}}(-X),\textbf{S}_{01}=\sigma_{\textbf{S}_{00},\textbf{S}_{11}}(-Y)$.
Since $\mbox{ker}\, \sigma_{\textbf{T}_{00},\textbf{S}_{00}}=\mbox{ker}\, \sigma_{\textbf{S}_{11},\textbf{T}_{11}}=\{0\}$ and $\widetilde{\textbf{T}}\sim_{u}\widetilde{\textbf{S}}$,
by Lemma \ref{mainlemma}, we will find a unitary operator $U=(\!(U_{i,j})\!)_{2\times 2}$ such that
\begin{equation}\label{21.10.22}
T_{00}^i=U_{10}^{-1}S_{11}^iU_{10}, T_{11}^i=U^*_{01}S_{00}^iU^{*-1}_{01},
\end{equation}
and $S_{00}^i(Y-U_{01}X^*U^{-1}_{10})=(Y-U_{01}X^*U^{-1}_{10})S_{11}^i,\,1\leq i\leq m.$
 Multiplying $U_{10}$ on the right side of the equation above, by equation (\ref{21.10.22}), we have
\begin{eqnarray}\label{mp14}
S_{00}^i(YU_{10}-U_{01}X^*)=(YU_{10}-U_{01}X^*)T_{00}^i,1\leq i\leq m.
\end{eqnarray}
Then multiplying $U^*_{01}$ on the left side of the last equation above, by equation (\ref{21.10.22}) again and $XU^*_{01}U_{01}=U^*_{10}U_{10}X$ due to Lemma \ref{mainlemma}, we obtain
\begin{eqnarray}\label{mp15}
T_{11}^i( U^*_{01}Y-X^*U^*_{10})=(U^*_{01}Y-X^*U^*_{10})S_{11}^i,\,1\leq i\leq m.
\end{eqnarray}

Set $S_0=YU_{10}-U_{01}X^*$ and $S_1=U^*_{01}Y-X^*U^*_{10}$. By equations (\ref{mp14}) and (\ref{mp15}), we see that $S_0\in \mbox{ker} \sigma_{\textbf{S}_{00},\textbf{T}_{00}}$ and
$S_1\in \mbox{ker} \sigma_{\textbf{T}_{11},\textbf{S}_{11}}$. Let $\widehat{T}_i=\left(\begin{smallmatrix}
 S_{00}^i & S_0 \\
 0 & T_{00}^i\\
 \end{smallmatrix}\right),\widehat{S}_i=\left(\begin{smallmatrix}
 T_{11}^i & S_1 \\
 0 & S_{11}^i\\
 \end{smallmatrix}\right)$, $1\leq i\leq m$ and $\widehat{\textbf{T}}=(\widehat{T}_1,\cdots,\widehat{T}_m),\widehat{\textbf{S}}=(\widehat{S}_1,\cdots,\widehat{S}_m)$.
That means $\widehat{\textbf{T}}, \widehat{\textbf{S}}\in \mathcal{FB}^{m}_2(\Omega)$ from Definition \ref{220418.1}.
Set $Z:=U^*_{01}\oplus U_{10}$. Then $Z$ is invertible. Using the equations $XU^*_{01}U_{01}=U^*_{10}U_{10}X$ and (\ref{21.10.22}) again, we imply that
$Z\widehat{T}_iZ^{-1}=\widehat{S}_i$
for all $1\leq i\leq m$. Hence, $Z\widehat{\textbf{T}}=\widehat{\textbf{S}}Z$ and $\widehat{\textbf{T}}\sim_s \widehat{\textbf{S}}$.
\end{proof}

Let $\{\textbf{T}\}^{\prime}=\{X|X\textbf{T}=\textbf{T}X\}$, $\{\textbf{T},\textbf{T}^*\}^{\prime}=\{X|X\textbf{T}=\textbf{T}X,\,X\textbf{T}^*=\textbf{T}^*X\}.$ The commuting tuple $\textbf{T}$ is said to be irreducible,  if there is no nontrivial orthogonal idempotents in $\{\textbf{T}\}^{\prime}$.
The following lemma is given by J. Fang, C. Jiang and P. Wu in Lemma 3.3 of \cite{FJW6}, which shows that the double commutant $\{T,T^*\}^{\prime}$ of irreducible operator $T$ is only scalar operators.
We will prove that this result also holds for irreducible operator tuples.

\begin{lem}\cite{FJW6}\label{21.9.20}
If $\text{\bf{T}}\in{\mathcal L}({\mathcal H})^{m}$ is irreducible and there is $X\in{\mathcal L}({\mathcal H})$ such that $X\in\{\text{\bf{T}},\text{\bf{T}}^*\}^{\prime}$,
then $X$ is a scalar multiple of identity.
\end{lem}
\begin{proof}
Since $X\textbf{T}=\textbf{T}X,\,X\textbf{T}^*=\textbf{T}^*X$, we have $X^*X\textbf{T}=\textbf{T}X^*X$. Then, for any spectral projection $P$ of $X^*X$,
this implies $P\textbf{T}=\textbf{T}P$. From the irreducibility of $\textbf{T}$, it follows that $P=0$ or $I$.
Furthermore, $\sigma(X^*X)=\{\alpha\}$ and $X^*X=\alpha I$.
Note that $X\textbf{T}(\mbox{ker}X)=\textbf{T}X(\mbox{ker}X)=0$ and $X\textbf{T}^*(\mbox{ker}X)=\textbf{T}^*X(\mbox{ker}X)=0$.
We know that $\mbox{ker}X$ is a reducing subspace for $\textbf{T}$, then $\mbox{ker}X=\{0\}$ or ${\mathcal H}$, since $\textbf{T}$ is irreducible.
So, either $X$ is injective or $X$ is 0. Suppose that $X$ is a injective in the dense range. By the polar decomposition of $X$, we have $X=U(X^*X)^{\frac{1}{2}}$,
$U$ is a unitary operator.
We assume that $\alpha\neq0$, then $U\textbf{T}=\textbf{T}U, U\textbf{T}^*=\textbf{T}^*U$. Repeating the above assumption, we
have $U=\beta I$. Thus $X=\sqrt{\alpha}\beta I$ is a scalar multiple of identity.
\end{proof}

Let $T_1,T_2$ be two bounded operators acting on $\mathcal{H}$ and $\alpha\in \mathbb{C}$.
For $\widetilde{T}=\left(\begin{smallmatrix}
 T_{1} & \alpha(T_{1}-T_2) \\
 0 & T_{2} \\
\end{smallmatrix}\right)$ on $\mathcal{H}\oplus\mathcal{H}$, in Lemma 2.1 of \cite{AK}, A. Koranyi proved that the operator $\tilde{T}$ is unitarily equivalent to $\left(\begin{smallmatrix}
 T_{2} & \alpha(T_{2}-T_1) \\
 0 & T_{1} \\
\end{smallmatrix}\right)$ through intertwining unitary operator $\frac{1}{\sqrt{1+\alpha^2}}\left(\begin{smallmatrix}
 -\alpha I & I \\
 I & \alpha I \\
\end{smallmatrix}\right)$. In the following two propositions, we discuss the conditions that the tuples in $N\mathcal{FB}_{n_0,n_1}$ make this conclusion hold, which is similar to the above.

\begin{prop}
For $i=0,1$, let $\text{\bf{T}}_{ii}\in \mathcal{B}_{n}^{m}(\Omega)\cap \mathcal{L}(\mathcal{H})^m$ and $\text{\bf{T}}_{01}=\sigma_{\text{\bf{T}}_{00},\text{\bf{T}}_{11}}(-X),\text{\bf{S}}_{01}=\sigma_{\text{\bf{T}}_{11},\text{\bf{T}}_{00}}(-Y)$ for $X,Y\in\mathcal{L}(\mathcal{H})$. Let
$\widetilde{\text{\bf{T}}}=(T_{1},\cdots, T_{m}),\widetilde{\text{\bf{S}}}=(S_{1},\cdots, S_{m})\in N\mathcal{FB}_{2n}^{m}(\Omega)$ with
$T_{i}=\left(\begin{smallmatrix}
 T_{00}^i & T_{01}^i \\
 0 & T_{11}^i \\
\end{smallmatrix}\right),\ S_{i}=\left(\begin{smallmatrix}
T_{11}^i & S_{01}^i \\
 0 & T_{00}^i \\
\end{smallmatrix}\right),\,1\leq i\leq m.$
Suppose that $\mbox{ker}\,\sigma_{\text{\bf{T}}_{00},\text{\bf{T}}_{11}}=\{0\}$, $\text{\bf{T}}_{00},\text{\bf{T}}_{11}$ are irreducible and $XX^*\in\{\text{\bf{T}}_{00}\}^{\prime}$, $X^*X\in\{\text{\bf{T}}_{11}\}^{\prime}$.
Then $\widetilde{\text{\bf{T}}}\sim_{u} \widetilde{\text{\bf{S}}}$ if and only if there exists $\theta\in \mathbb{R},$ such that
$\text{\bf{S}}_{01}=e^{i\theta}\ \sigma_{\text{\bf{T}}_{11},\text{\bf{T}}_{00}}(-X^*).$
\end{prop}

\begin{proof}
Let $U=(\!(U_{i,j})\!)_{2\times 2}$ be a unitary operator which satisfies that $U\widetilde{\textbf{T}}=\widetilde{\textbf{S}}U$.
By Lemma \ref{mainlemma}, we have $U_{01},\,U_{10}$ are invertible and $(I+X^*X)^{-1}=U^*_{01}U_{01}$,
$(I+XX^*)^{-1}=U^*_{10}U_{10}$. Since $X$ is a bounded linear operator, then $I+X^*X$ and $I+XX^*$ are positive and invertible.
Furthermore, we have that $U_{1}:=(I+X^*X)^{\frac{1}{2}}U_{01}^*,\,U_{2}:=U_{10}(I+XX^*)^{\frac{1}{2}}$ are unitary.
By using the statement (1) of Lemma \ref{mainlemma}, we also have
$U_{10}T_{00}^i=T_{00}^iU_{10},\,T_{11}^iU_{01}^*=U_{01}^*T_{11}^i,\,1\leq i\leq m.$ It follows that
\begin{equation}\label{21.910}
U_{2}(I+XX^*)^{-\frac{1}{2}}T_{00}^i=T_{00}^iU_{2}(I+XX^*)^{-\frac{1}{2}}
\end{equation}
and
\begin{equation}\label{21.9101}
T_{11}^i(I+X^*X)^{-\frac{1}{2}}U_{1}=(I+X^*X)^{-\frac{1}{2}}U_{1}T_{11}^i,\,1\leq i\leq m.
\end{equation}
From the conditions $XX^*\in\{T_{00}^i\}^{\prime}$, $X^*X\in\{T_{11}^i\}^{\prime}$, we obtain that
$(I+XX^*)T_{00}^i=T_{00}^i(I+XX^*)$ and $(I+X^*X)T_{11}^i=T_{11}^i(I+X^*X),\,1\leq i\leq m.$ By functional calculus of positive operators, we have
$$(I+XX^*)^{-\frac{1}{2}}T_{00}^i=T_{00}^i(I+XX^*)^{-\frac{1}{2}},\,(I+X^*X)^{-\frac{1}{2}}T_{11}^i=T_{11}^i(I+X^*X)^{-\frac{1}{2}},\,1\leq i\leq m.$$
Combining with equations (\ref{21.910}) and (\ref{21.9101}), we imply that $U_{2}T_{00}^i=T_{00}^iU_{2},\,U_{1}T_{11}^i=T_{11}^iU_{1},\,1\leq i\leq m.$ Thus,
$U_{2}\in\{\textbf{T}_{00},\textbf{T}_{00}^*\}^{\prime},\,U_{1}\in\{\textbf{T}_{11},\textbf{T}_{11}^*\}^{\prime}$, since $U_{1},U_{2}$ are unitary.
From $\textbf{T}_{00},\textbf{T}_{11}$ are irreducible and Lemma \ref{21.9.20}, we obtain
$U_{1}=e^{i\theta_1}I, U_{2}=e^{i\theta_2}I$ for some $\theta_1,\theta_2\in \mathbb{R}.$
By the statement (3) of Lemma \ref{mainlemma}, we have $Y-U_{01}X^*U^{-1}_{10}\in \mbox{ker}\ {\sigma}_{\textbf{T}_{11},\textbf{T}_{00}}$. It follows that
$T_{11}^i(Y-U_{01}X^*U^{-1}_{10})=T_{11}^i(Y-e^{-i(\theta_1+\theta_2)}X^*)=(Y-e^{-i(\theta_1+\theta_2)}X^*)T_{00}^i$
and $YT_{00}^i-T_{11}^iY=e^{-i(\theta_1+\theta_2)}(X^*T_{00}^i-T_{11}^iX^*),\,1\leq i\leq m.$ This finishes the proof of necessary part.

For the proof of the sufficient part, choose any $
\theta_1, \theta_2 \in \mathbb{R}$ such that $\theta_1+\theta_2=-\theta$. Define the operator $U$ as follows
$$U=\begin{pmatrix} e^{-i\theta_1}(I+X^*X)^{-\frac{1}{2}}X^*&e^{-i\theta_1}(I+X^*X)^{-\frac{1}{2}}\\ e^{i\theta_2}(I+XX^*)^{-\frac{1}{2}}&-e^{i\theta_2}(I+XX^*)^{-\frac{1}{2}}X\end{pmatrix}.$$
Using the fact of $X(I+X^*X)^{-1}=(I+XX^*)^{-1}X$, we obtain $U$ is a unitary operator. By a simple calculation, we imply $UT_{i}=S_{i}U$ for $1\leq i\leq m$,
then $U\widetilde{\textbf{T}}=\widetilde{\textbf{S}}U$.
\end{proof}

\begin{prop}
For $i=0,1$, let $\text{\bf{T}}_{ii}\in \mathcal{B}_{n}^{m}(\Omega)$ and $\text{\bf{T}}_{01}=\sigma_{\text{\bf{T}}_{00},\text{\bf{T}}_{11}}(-X)$ for some self-adjoint operator $X$. Let
$\widetilde{\text{\bf{T}}}=(T_{1},\cdots, T_{m})\in N\mathcal{FB}_{2n}^{m}(\Omega)$ with
$T_{i}=\left(\begin{smallmatrix}
 T_{00}^i & T_{01}^i \\
 0 & T_{11}^i \\
\end{smallmatrix}\right),\,1\leq i\leq m.$
Suppose that $X\in\{\text{\bf{T}}_{00}\}^{\prime}\cap\{\text{\bf{T}}_{11}\}^{\prime}$.
Then the operator $\widehat{\text{\bf{T}}}$ gotten by interchanging the roles of $\text{\bf{T}}_{00}$ and $\text{\bf{T}}_{11}$ is unitarily equivalent to $\widetilde{\text{\bf{T}}}$.
\end{prop}
\begin{proof}
Let $\widetilde{X}=\left(\begin{smallmatrix}
 X(I+X^2)^{-\frac{1}{2}} & (I+X^2)^{-\frac{1}{2}} \\
 (I+X^2)^{-\frac{1}{2}} & -X(I+X^2)^{-\frac{1}{2}} \\
\end{smallmatrix}\right)$. By functional calculus of positive operators, we have $X(I+X^2)^{-\frac{1}{2}}=(I+X^2)^{-\frac{1}{2}}X$, thus $\widetilde{X}$ is self-adjoint. Note that
$$\widetilde{X}\widetilde{X}^*=\widetilde{X}^*\widetilde{X}=\left(\begin{smallmatrix}
 X(I+X^2)^{-1}X+(I+X^2)^{-1} & X(I+X^2)^{-1}-(I+X^2)^{-1}X \\
 (I+X^2)^{-1}X-X(I+X^2)^{-1} & (I+X^2)^{-1}+X(I+X^2)^{-1}X \\
\end{smallmatrix}\right)=\left(\begin{smallmatrix}
 I & 0 \\
 0 & I \\
\end{smallmatrix}\right),$$
we see that $\widetilde{X}$ is unitary. From $X\in\{\textbf{T}_{00}\}^{\prime}\cap\{\textbf{T}_{11}\}^{\prime}$, then for any $1\leq i\leq m$, we have
$XT_{jj}^i=T_{jj}^iX,j=0,1$. By functional calculus of $I+X^2$, we also have $$(I+X^2)^{-\frac{1}{2}}T_{jj}^i=T_{jj}^i(I+X^2)^{-\frac{1}{2}}\,\ \mbox{and}\,\
(I+X^2)^{\frac{1}{2}}T_{jj}^i=T_{jj}^i(I+X^2)^{\frac{1}{2}},1\leq i\leq m,j=0,1. $$
Based on a simple calculation, $\widetilde{X}\widetilde{\textbf{T}}=\widehat{\textbf{T}}\widetilde{X}$ can be obtained. Hence, $\widetilde{\textbf{T}}$ is unitarily equivalent to $\widehat{\textbf{T}}$.
\end{proof}

\subsection{Some properties of tuples in $N\mathcal{F}B^m_{n_0,n_1}(\Omega)$}

The commuting operator tuple $\textbf{T}$ is said to be strongly irreducible if there is no nontrivial idempotents in $\{\textbf{T}\}^{\prime}$. Otherwise, it is strongly reducible.
A strongly irreducible operator can be regarded as a natural generalization of
a Jordan block matrix on the infinite dimensional case. In \cite{Jiang1}, C. Jiang proved that for any strongly irreducible Cowen-Douglas operator
$T$, $\{T\}^{\prime}/rad(\{T\}^{\prime})$ is commutative, where $rad(\{T\}^{\prime})$ denotes the Jacobson radical of
$\{T\}^{\prime}$. Based on this, C. Jiang gave a similarity classification of strongly irreducible Cowen-Douglas operators by
using the $K_0$-group of their commutant algebra as an invariant (see more details in \cite{Jiang1}). These results are also
generalized to the case of direct integrals of strongly irreducible operators by R. Shi (cf. \cite{Rui}).
The following proposition shows that the strong reducibility of tuples in $N\mathcal{FB}^{m}_{n_0,n_1}(\Omega)$, that is,
every tuple in $N\mathcal{FB}^{m}_{n_0,n_1}(\Omega)$ can be written as the direct sum of two tuples in $\mathcal{B}^{m}_{n_i}(\Omega),i=0,1$ up to similarity.
For $m$-tuples $\textbf{T}_{00}$ and $\textbf{T}_{11}$, $\textbf{T}_{00}\oplus \textbf{T}_{11}=(T_{00}^1\oplus T_{11}^1,\cdots,T_{00}^m\oplus T_{11}^m)$.
\begin{prop}\label{3.5pro}
For $i=0,1$, let $\text{\bf{T}}_{ii}\in \mathcal{B}_{n_i}^{m}(\Omega)\cap \mathcal{L}(\mathcal{H}_i)^m$ and $\text{\bf{T}}_{01}=\sigma_{\text{\bf{T}}_{00},\text{\bf{T}}_{01}}(-X)$ for $X\in\mathcal{L}(\mathcal{H}_1,\mathcal{H}_0)$. Let
$\widetilde{\text{\bf{T}}}=(T_{1},\cdots, T_{m})\in N\mathcal{FB}_{n_0,n_1}^{m}(\Omega)$ with
$T_{i}=\left(\begin{smallmatrix}
 T_{00}^i & T_{01}^i \\
 0 & T_{11}^i \\
\end{smallmatrix}\right),\,1\leq i\leq m.$
Then $\widetilde{\text{\bf{T}}}$ is strongly reducible. What is more, $\widetilde{\text{\bf{T}}}$ is similar to $\text{\bf{T}}_{00}\oplus \text{\bf{T}}_{11}$.
\end{prop}

\begin{proof}

Let $W=\left(\begin{smallmatrix}I&-X\\0&I\end{smallmatrix}\right)$. We have that
$$WT_{j}=\begin{pmatrix}T_{00}^j&-T_{00}^jX\\0&T_{11}^j\end{pmatrix}
=\begin{pmatrix}T_{00}^j&0\\0&T_{11}^j\end{pmatrix}\begin{pmatrix}I&-X\\0&I\end{pmatrix}=(T_{00}^j\oplus T_{11}^j)W,\,1\leq j\leq m $$
and $W\widetilde{\textbf{T}}=(\textbf{T}_0\oplus \textbf{T}_1)W$.
Note that $W$ is invertible and $W^{-1}=\left(\begin{smallmatrix}I&X\\0&I\end{smallmatrix}\right).$ Then we finish the proof.
\end{proof}
The characterization  of irreducibility of tuples in $N\mathcal{FB}^{m}_{2}(\Omega)$ is as follows.
\begin{prop}
Let $\text{\bf{T}}_{ii}\in\mathcal{B}^{m}_1(\Omega)$, $\text{\bf{T}}_{ii}\sim_{u}(\text{\bf{M}}^*_z, \mathcal{H}_{K_{i}}),i=0,1$ and $\text{\bf{T}}_{01}=\sigma_{\text{\bf{T}}_{00},\text{\bf{T}}_{11}}(-X)$ for some $X$.
Suppose that
$\widetilde{\text{\bf{T}}}=(T_{1}, \cdots, T_{m})\in N\mathcal{FB}_{2}^{m}(\Omega)$ with $T_{i}=\left(\begin{smallmatrix}
 T_{00}^i & T_{01}^i \\
 0 & T_{11}^i \\
 \end{smallmatrix}\right),1\leq i\leq m.$
If $\lim\limits_{dist(w, \partial \Omega)\rightarrow 0}\frac{K_0(w,w)}{K_1(w,w)}=0,$
then $\widetilde{\text{\bf{T}}}$ is irreducible.
\end{prop}
\begin{proof}
Suppose that $\widetilde{\textbf{T}}$ is reducible, then there exists a nontrivial orthogonal projection $P=\left (\begin{smallmatrix}
 P_{00} & P_{01} \\
 P_{10} & P_{11} \\
\end{smallmatrix}\right)\in\{\widetilde{\textbf{T}}\}^{\prime}$, such that
\begin{equation}\label{21715.1}
\left (\begin{smallmatrix}
 P_{00}T_{00}^i & P_{00}(XT_{11}^i-T_{00}^iX)+P_{01}T_{11}^i \\
 P_{10}T_{00}^i & P_{10}(XT_{11}^i-T_{00}^iX)+P_{11}T_{11}^i \\
\end{smallmatrix}\right )=\left (\begin{smallmatrix}
 T_{00}^iP_{00}+(XT_{11}^i-T_{00}^iX)P_{10} & T_{00}^iP_{01}+(XT_{11}^i-T_{00}^iX)P_{11} \\
 T_{11}^iP_{10} & T_{11}^iP_{11} \\
\end{smallmatrix}\right )
\end{equation}
for all $1\leq i\leq m$. It follows that $P_{10}\in \bigcap_{i=1}^{m}\mbox{ker}\ \sigma_{T_{11}^i,T_{00}^i}=\mbox{ker}\ \sigma_{\textbf{T}_{11},\textbf{T}_{00}}$.

By Lemma \ref{intertwining}, if we have $\lim\limits_{dist(w, \partial \Omega)\rightarrow 0}\frac{K_0(w,w)}{K_1(w,w)}=0,$ then $\mbox{ker}\ \sigma_{\textbf{T}_{11},\textbf{T}_{00}}=\{0\}$ and $P_{10}=0$.
Note that $P$ is a self-adjoint idempotent, we obtain $P_{01}=0$ and $P_{ii}=P_{ii}^*=P_{ii}^2,\,i=0,1$.
From equation (\ref{21715.1}), we infer
$P_{00}T_{00}^i=T_{00}^iP_{00},\,P_{11}T_{11}^i=T_{11}^iP_{11},\,1\leq i\leq m.$
Then $P_{00}\textbf{T}_{00}=\textbf{T}_{00}P_{00},\,P_{11}\textbf{T}_{11}=\textbf{T}_{11}P_{11}$.
Since tuples in $\mathcal{B}_1^m(\Omega)$ are irreducible, we have $P_{ii}=0$ or $I$.
According to $P_{00}(XT_{11}^i-T_{00}^iX)=(XT_{11}^i-T_{00}^iX)P_{11},$ we have $P_{00}=P_{11}=0$ or $I$, that is, $P$ is trivial. This is a contradiction.
Hence, $\widetilde{\textbf{T}}$ is irreducible.
\end{proof}
By the following proposition, the Hermitian holomorphic vector bundles corresponding to the tuples in $N\mathcal{FB}_{n_0,n_1}^{m}(\Omega)$ is given.
\begin{prop}\label{frame}
Let $\text{\bf{T}}_{ii}\in\mathcal{B}^{m}_{n_i}(\Omega)$, $\text{\bf{T}}_{ii}\sim_{u}(\text{\bf{M}}^*_z, \mathcal{H}_{K_{i}}),i=0,1$ and $\text{\bf{T}}_{01}=\sigma_{\text{\bf{T}}_{00},\text{\bf{T}}_{11}}(-X)$ for some $X$.
Suppose that
$\widetilde{\text{\bf{T}}}=(T_{1}, \cdots, T_{m})\in N\mathcal{FB}_{n_0+n_1}^{m}(\Omega)$ with $T_{i}=\left(\begin{smallmatrix}
 T_{00}^i & T_{01}^i \\
 0 & T_{11}^i \\
 \end{smallmatrix}\right),1\leq i\leq m.$
Then for all $w\in \Omega$,
$$E_{\text{\bf{T}}}(w)=\overline{span}\{K_0(\cdot,\bar{w})\xi_0,X(K_1(\cdot,\bar{w})\xi_1)+K_1(\cdot,\bar{w})\xi_1,\xi_0\in \mathbb{C}^{n_0},\xi_1\in \mathbb{C}^{n_1}\}.$$
\end{prop}
\begin{proof}
Since $E_{\textbf{T}_{ii}}(w)=\overline{span}\{K_i(\cdot,\bar{w})\xi_i,\xi_i\in \mathbb{C}^{n_i}\}$ and the dimension of $E_{\mathbf{T}_{ii}}(w)$ is $n_i$, $i=0,1$, it is easy to see that $$K_0(\cdot,\bar{w})\xi_0, X(K_1(\cdot,\bar{w})\xi_1)+K_1(\cdot,\bar{w})\xi_1\in\ker(\widetilde{\textbf{T}}-w),w\in \Omega, \xi_0\in \mathbb{C}^{n_0},\xi_1\in \mathbb{C}^{n_1}. $$
Note that dim$\,\ker(\widetilde{\textbf{T}}-w)=n_0+n_1,w\in\Omega$,
then we only need to prove that for each $\xi_0\in \mathbb{C}^{n_0},\xi_1\in \mathbb{C}^{n_1}$, $K_0(\cdot,\bar{w})\xi_0$ and $X(K_1(\cdot,\bar{w})\xi_1)+K_1(\cdot,\bar{w})\xi_1$ are linearly independent.
For fixed but arbitrary $\xi_0\in \mathbb{C}^{n_0},\xi_1\in \mathbb{C}^{n_1}$, suppose that there exist $x_0,x_1\in \mathbb{C}$ such that $$x_0K_0(\cdot,\bar{w})\xi_0+x_1(X(K_1(\cdot,\bar{w})\xi_1)+K_1(\cdot,\bar{w})\xi_1)=0.$$
By taking the inner product with $K_1(\cdot,\bar{w})\xi',\xi'\in\mathbb{C}^{n_1}$ on both sides, we have that
$\big\langle x_1K_1(\cdot,\bar{w})\xi_1, K_1(\cdot,\bar{w})\xi'\big\rangle=0$.
With the spanning property of $\{K_1(\cdot,\bar{w})\xi',\xi'\in\mathbb{C}^{n_1}\}$, we infer $x_1K_1(\cdot,\bar{w})\xi_1=0$, then $x_{1}=0$, since $K_1(\cdot,\bar{w})\xi_1$ is non-zero.
Thus, we obtain $x_0K_0(\cdot,\bar{w})\xi_0=0$. Similarly, we have $x_0=0$. This completes the proof.
\end{proof}

\begin{ex}
For $i=0,1$, let $\textbf{T}_{ii},\textbf{S}_{ii}\in\mathcal{B}_{1}^{m}(\Omega)$, $\textbf{T}_{ii}=(\mathbf{M}_z^*,\mathcal{H}_{K_i})$, $\textbf{S}_{ii}=(\mathbf{M}_z^*,\mathcal{H}_{\widetilde{K}_i})$ and $\textbf{T}_{01}=\sigma_{\textbf{T}_{00},\textbf{T}_{11}}(-X),\textbf{S}_{01}=\sigma_{\textbf{S}_{00},\textbf{S}_{11}}(-Y)$ for some $X,Y$. Let
$\widetilde{\textbf{T}}=(T_{1},\cdots, T_{m}),\widetilde{\textbf{S}}=(S_{1}, \cdots, S_{m})\in N\mathcal{FB}_{2}^{m}(\Omega)$ with
$T_i=\left(\begin{smallmatrix}T_{00}^i&T_{01}^i\\0&T_{11}^i\end{smallmatrix}\right)$, $S_i=\left(\begin{smallmatrix}S_{00}^i&S_{01}^i\\0&S_{11}^i\end{smallmatrix}\right), 1\leq i\leq m$.
By Lemma \ref{frame}, we have $\{K_0(\cdot,\bar{w}), X(K_1(\cdot,\bar{w}))+K_1(\cdot,\bar{w})\}$ is a frame of $E_{\widetilde{\textbf{T}}}(w)$. Similarly, a frame of $E_{\widetilde{\textbf{S}}}(w)$ is obtained.

Define
$K_{\gamma}, K_{\widetilde{\gamma}}$ to be the function on $\Omega^*\times \Omega^*$
taking values in the $2\times 2$ matrices
$\mathcal{M}_2(\mathbb{C})$:
$$K_{\gamma}(z,w)= \begin{pmatrix}
  K_0(z,w) &
  \langle X(K_1(\cdot,w)),K_0(\cdot, z)\rangle \\
 \langle K_0(\cdot, w), X(K_1(\cdot, z))\rangle
  &
 \langle X(K_1(\cdot, w)), X(K_1(\cdot, z))\rangle +K_1(z,w)\\
\end{pmatrix},$$
$$K_{\widetilde{\gamma}}(z,w)=\begin{pmatrix}
  \widetilde{K}_0(z,w) &
  \langle Y(\widetilde{K}_1(\cdot, w)),\widetilde{K}_0(\cdot, z)\rangle \\
 \langle \widetilde{K}_0(\cdot, w), Y(\widetilde{K}_1(\cdot, z))\rangle
  &
 \langle Y(\widetilde{K}_1(\cdot, w)), Y(\widetilde{K}_1(\cdot, z))\rangle +\widetilde{K}_1(z,w)\\
\end{pmatrix}.$$

By subsection 2.2 in \cite{Dinesh}, we know that $\widetilde{\textbf{T}}$ and $\widetilde{\textbf{S}}$ are unitarily equivalent to the adjoint of multiplication operator tuple $\textbf{M}_z$ on some analytic functional spaces $\mathcal{H}_{K_{\gamma}}$ and  $\mathcal{H}_{K_{\widetilde{\gamma}}}$ with
reproducing kernel $K_{\gamma}(z,w)$ and $K_{\widetilde{\gamma}}(z,w)$, respectively. That means $\widetilde{\textbf{T}}\sim_{u} (\textbf{M}^*_z, \mathcal{H}_{K_{\gamma}}), \widetilde{\textbf{S}}\sim_{u} (\textbf{M}^*_z, \mathcal{H}_{\widetilde{K}_{\gamma}})$.
R.E. Curto and N. Salinas gave a necessary and sufficient condition for the unitary equivalence of commuting operator tuples acting on reproducing kernel Hilbert spaces (see Remark 3.8, \cite{CS2}),
that is,
$\textbf{M}_z$ acting on $\mathcal{H}_{K_{\gamma}}$ and $\mathcal{H}_{K_{\widetilde{\gamma}}}$ are unitarily equivalent if and only if $\Phi(z)K_{\gamma}(z,w)\overline{\Phi^{T}(w)}=K_{\widetilde{\gamma}}(z,w)$ for some holomorphic and invertible function $\Phi$.

Now if there exist holomorphic functions $\phi$ and $\psi$ such that $\Phi(w):=\left (\begin{smallmatrix}0&\phi(w)\\
\psi(w)&0
\end{smallmatrix}\right)$ which satisfies $\Phi(z)K_{\gamma}(z,w)\overline{\Phi^{T}(w)}=K_{\widetilde{\gamma}}(z,w)$,
then $\textbf{T}$ is unitarily equivalent to $\widetilde{\textbf{T}}$, that is,
\begin{eqnarray*}
&&\left(\begin{smallmatrix}0&\phi(z)\\
\psi(z)&0
\end{smallmatrix}\right)
\left(\begin{smallmatrix}
  K_0(z,w) &
  \langle X(K_1(\cdot, w)),K_0(\cdot, z)\rangle \\
 \langle K_0(\cdot, w), X(K_1(\cdot, z))\rangle
  &
 \langle X(K_1(\cdot, w)), X(K_1(\cdot, z))\rangle +K_1(z,w)\\
\end{smallmatrix}\right)\left(\begin{smallmatrix}0&\overline{\psi(w)}\\
\overline{\phi(w)}&0
\end{smallmatrix}\right)\\
&=&\left(\begin{smallmatrix}
\phi(z)(\langle X(K_1(\cdot, w)), X(K_1(\cdot, z))\rangle + K_1(z, w))\overline{\phi(w)}  &
\phi(z)(\langle K_0(\cdot, w), X(K_1(\cdot, z))\rangle) \overline{\psi(w)} \\ \psi(z)( \langle X(K_1(\cdot, w)),K_0(\cdot, z)\rangle)\overline{\phi(w)} & \psi(z)K_0(z, w)\overline{\psi(w)}
 \\
\end{smallmatrix}\right)\\
&=& \left(\begin{smallmatrix}
  \widetilde{K}_0(z,w) &
  \langle Y(\widetilde{K}_1(\cdot, w)),\widetilde{K}_0(\cdot, z)\rangle \\
 \langle \widetilde{K}_0(\cdot, w), Y(\widetilde{K}_1(\cdot, z))\rangle
  &
 \langle Y(\widetilde{K}_1(\cdot, w)), Y(\widetilde{K}_1(\cdot, z))\rangle +\widetilde{K}_1(z,w)\\
\end{smallmatrix}\right).\\
\end{eqnarray*}
Choosing $z=w$, we have that $$\widetilde{K}_0(w,w)=\|X(\phi(w)K_1(\cdot,w))\|^2+|\phi(w)|^2K_1(w,w)=\|(I+X^*X)^{\frac{1}{2}}(\phi(w)K_1(\cdot,w))\|^2, $$
$$|\psi(w)|^2K_0(w,w)=\|Y(\widetilde{K}_1(\cdot,w))\|^2+\widetilde{K}_1(w,w)=\|(I+Y^*Y)^{\frac{1}{2}}\widetilde{K}_1(\cdot,w)\|^2,w\in \Omega.$$
By Theorem \ref{main9.26} and the proof of Example \ref{ALex},
we have that $\textbf{T}_{00}\sim_{s} \textbf{S}_{11}, \textbf{S}_{00}\sim_{s} \textbf{T}_{11}$. Thus, $\textbf{T}_{00}\oplus \textbf{S}_{00} \sim_{s} \textbf{T}_{11}\oplus \textbf{S}_{11}.$
\end{ex}
Let $T,S\in \mathcal{B}^1_1(\Omega)$ and $T\sim_{u} (M^*_z, \mathcal{H}_{K_{0}}),S\sim_{u} (M^*_z, \mathcal{H}_{K_{1}})$.
By Lemma \ref{intertwining}, we know that
if $T_0\sim_{s}T_1$, then $\frac{K_0(w,w)}{K_1(w,w)}$ is bounded and and bounded below from zero.
In the following proposition, we will prove that there is a similar result in the operator class $N\mathcal{FB}^1_{2}(\Omega)$.
For the case of index two, $\frac{K_0(w,w)}{K_1(w,w)}$ is replaced by the ratio of the determinants of the metrics corresponding to the two bundles.

\begin{prop}
Let $T_{ii},S_{ii}\in \mathcal{B}^1_1(\Omega)$ and $T_{ii}\sim_{u} (M^*_z, \mathcal{H}_{K_{i}}), S_{ii}\sim_{u} (M^*_z, \mathcal{H}_{\widetilde{K}_{i}}), i=0,1$.
Suppose that $T=\left(\begin{smallmatrix}
 T_{00} & T_{01} \\
 0 & T_{11} \\
 \end{smallmatrix}\right),\ S=\left(\begin{smallmatrix}
 S_{00} & S_{01} \\
 0 & S_{11} \\
 \end{smallmatrix}\right)\in N\mathcal{FB}^1_{2}(\Omega)$ and there exist $X,Y$ such that $T_{01}=\sigma_{T_{00},T_{11}}(-X),S_{01}=\sigma_{S_{00},S_{11}}(-Y)$.
 If
$T\sim_{s}S$, then there exist metrics $h_{T},h_{S}$ corresponding to $E_T, E_{S}$ such that
$m\leq \frac{det\ h_{T}(\omega)}{det\ h_{S}(\omega)}\leq M, \omega\in\Omega,$ for positive numbers $m$ and $M$.
\end{prop}

\begin{proof}
Without loss of generality, we assume that $T_{ii}= (M^*_z, \mathcal{H}_{K_{i}}), S_{ii}= (M^*_z, \mathcal{H}_{\widetilde{K}_{i}}), i=0,1$.
Then $K_i(\cdot,\bar{w}), \widetilde{K}_i(\cdot,\bar{w})$ are the sections of $E_{T_{ii}}$ and $E_{S_{ii}},i=0,1$, respectively.
By Lemma \ref{frame}, we know that $\{K_0(\cdot,\bar{w}), XK_1(\cdot,\bar{w})+K_1(\cdot,\bar{w})\},\ \{\widetilde{K}_0(\cdot,\bar{w}), Y\widetilde{K}_1(\cdot,\bar{w})+\widetilde{K}_1(\cdot,\bar{w})\}$ are frames of
$E_{T}(w),E_{S}(w)$, respectively.
It follows that
$$h_{T}(\omega)=\left(\begin{matrix}K_0(\bar{w},\bar{w}) & \langle X(K_1(\cdot,\bar{w})),K_0(\cdot,\bar{w})\rangle\\
\langle K_0(\cdot,\bar{w}),X(K_1(\cdot,\bar{w}))\rangle & \|X(K_1(\cdot,\bar{w}))\|^{2}+K_1(\bar{w},\bar{w})\\
\end{matrix}\right ),$$
$$det\ h_{T}(\omega)=K_0(\bar{w},\bar{w})(K_1(\bar{w},\bar{w})+\|X(K_1(\cdot,\bar{w}))\|^{2})-|\langle K_0(\cdot,\bar{w}),X(K_1(\cdot,\bar{w}))\rangle|^{2}.$$ Similarly, we have
$$det\ h_{S}(\omega)=\widetilde{K}_0(\bar{w},\bar{w})(\widetilde{K}_1(\bar{w},\bar{w})+ \|Y(\widetilde{K}_1(\cdot,\bar{w}))\|^{2})-|\langle \widetilde{K}_0(\cdot,\bar{w}),Y(\widetilde{K}_1(\cdot,\bar{w}))\rangle|^{2}.$$

By Proposition \ref{3.5pro}, we know that operators in $N\mathcal{FB}^1_{2}(\Omega)$ are strongly reducible
and $T\sim_{s}T_{00}\oplus T_{11}$, $S\sim_{s}S_{00}\oplus S_{11}$. If $T\sim_{s}S$, then $T_{00}\oplus T_{11}\sim_{s}S_{00}\oplus S_{11}$.
By the main theorem of \cite{JGJ}, we know that every Cowen-Douglas operator has a unique strongly irreducible decomposition up to similarity. Thus, the equivalence relation is either
$T_{00}\sim_{s}S_{00}, T_{11}\sim_{s}S_{11}$ or $T_{00}\sim_{s}S_{11}, T_{11}\sim_{s}S_{00}$.
In either case, according to Lemma \ref{intertwining}, there exist positive numbers $m_{1}$ and $M_{1}$ such that
$m_{1}\leq\frac{det\ h_{T_{00}\oplus T_{11}}}{det\ h_{S_{00}\oplus S_{11}}}=\frac{K_0(\bar{w},\bar{w})K_1(\bar{w},\bar{w})}{\widetilde{K}_0(\bar{w},\bar{w})\widetilde{K}_1(\bar{w},\bar{w})}\leq M_{1}$.
By using Cauchy-Schwarz inequality, we have
$K_0(\bar{w},\bar{w})\|X(K_1(\cdot,\bar{w}))\|^{2}-|\langle K_0(\cdot,\bar{w}),X(K_1(\cdot,\bar{w}))\rangle|^{2} \geq0$. Thus,
$$
\frac{det\ h_{T}(\omega)}{det\ h_{T_{00}\oplus T_{11}}(\omega)}
=\frac{K_0(\bar{w},\bar{w})(K_1(\bar{w},\bar{w})+\|X(K_1(\cdot,\bar{w}))\|^{2})-|\langle K_0(\cdot,\bar{w}),X(K_1(\cdot,\bar{w}))\rangle|^{2}}{K_0(\bar{w},\bar{w})K_1(\bar{w},\bar{w})}
\geq1.$$
On the other hand, since $X$ is a bounded linear operator, we have
$$
\frac{det\ h_{T}(\omega)}{det\ h_{T_{00}\oplus T_{11}}(\omega)}
\leq\frac{K_0(\bar{w},\bar{w})(K_1(\bar{w},\bar{w})+\|X(K_1(\cdot,\bar{w}))\|^{2})}{K_0(\bar{w},\bar{w})K_1(\bar{w},\bar{w})}
\leq1+\|X\|^{2}.$$
Similarly, we obtain $1\leq\frac{det\ h_{S}}{det\ h_{S_{00}\oplus S_{11}}}\leq1+\|Y\|^{2}$.
Note that
$$\frac{det\ h_{T}}{det\ h_{S}}=\frac{det\ h_{T}}{det\ h_{T_{00}\oplus T_{11}}}\cdot \frac{det\ h_{T_{00}\oplus T_{11}}}{det\ h_{S_{00}\oplus S_{11}}}\cdot
\frac{det\ h_{S_{00}\oplus S_{11}}}{det\ h_{S}}.$$
Let $m:=\frac{m_{1}}{1+\|Y\|^{2}}, M:=M_{1}(1+\|X\|^{2})$. We have
$m\leq\frac{det\ h_{T}(\omega)}{det\ h_{S}(\omega)}\leq M.$ This completes the proof.
\end{proof}


\section{Weakly homogeneous operator tuples }\label{chapter5}
An
operator $T\in \mathcal{L}(\mathcal{H})$ is said to be weakly homogeneous if $\sigma(T)\subset\overline{\mathbb{D}}$ and $\phi(T)$ is similar to $T$
for each $\phi$ in $M\ddot{o}b$.
When $\mathcal{D}$ is a bounded symmetric domain, a commuting $m$-tuple
$\textbf{T}=(T_1,\cdots,T_m)$ of bounded operators is said to be weakly homogeneous with respect
to $G$ if
their joint Taylor spectrum is contained in $\overline{\mathcal{D}}$ and $\phi(\textbf{T})$ is similar to $\textbf{T}$ for every holomorphic automorphism $\phi\in G$.
Given a Hilbert space $\mathcal{H}$ with sharp reproducing kernel $K$ on $\mathbb{D}\times\mathbb{D}$, S. Ghara in \cite{Ghara} obtain an equivalent condition that the multiplication operator $M_z$ on $(\mathcal{H},K)$ is weakly homogeneous.
Next, we consider the weakly homogeneity of class $N\mathcal{FB}_{n_0+n_1}^{m}(\mathbb{D}^m)$.

\begin{prop}
Let $\text{\bf{T}}_{ii}\in\mathcal{B}^{m}_{n_i}(\mathbb{D}^m),i=0,1$ and $\text{\bf{T}}_{01}\in ran\sigma_{\text{\bf{T}}_{00},\text{\bf{T}}_{11}}$.
Suppose that
$\widetilde{\text{\bf{T}}}=(T_{1}, \cdots, T_{m})\in N\mathcal{FB}_{n_0+n_1}^{m}(\mathbb{D}^m)$ with $T_{i}=\left(\begin{smallmatrix}
 T_{00}^j & T_{01}^j \\
 0 & T_{11}^j \\
 \end{smallmatrix}\right),1\leq j\leq m.$
If $\text{\bf{T}}_{00},\text{\bf{T}}_{11}$ are both weakly homogeneous with respect to $M\ddot{o}b^m$, then $\widetilde{\text{\bf{T}}}$ is also $M\ddot{o}b^m$-weakly homogeneous.
\end{prop}
\begin{proof}
Suppose that there exists $X$ such that $\textbf{T}_{01}=\sigma_{\textbf{T}_{00},\textbf{T}_{11}}(-X)$.
By Proposition \ref{3.5pro}, we know that $\widetilde{\textbf{T}}$ is similar to $\textbf{T}_{00}\oplus \textbf{T}_{11}$ and $\left(\begin{smallmatrix}
 I & -X \\
 0 & I \\
 \end{smallmatrix}\right)\left(\begin{smallmatrix}
 T_{00}^j & T_{01}^j \\
 0 & T_{11}^j \\
 \end{smallmatrix}\right)\left(\begin{smallmatrix}
 I & X \\
 0 & I \\
 \end{smallmatrix}\right)=\left(\begin{smallmatrix}
 T_{00}^j & 0 \\
 0 & T_{11}^j \\
 \end{smallmatrix}\right)$ for $1\leq j\leq m$.
 Further, we have $\left(\begin{smallmatrix}
 I & -X \\
 0 & I \\
 \end{smallmatrix}\right)\left(\begin{smallmatrix}
 T_{00}^j & T_{01}^j \\
 0 & T_{11}^j \\
 \end{smallmatrix}\right)^n\left(\begin{smallmatrix}
 I & X \\
 0 & I \\
 \end{smallmatrix}\right)=\left(\begin{smallmatrix}
 T_{00}^j & 0 \\
 0 & T_{11}^j \\
 \end{smallmatrix}\right)^n$ for any positive integer $n$ and
 $$\left(\begin{smallmatrix}
 I & -X \\
 0 & I \\
 \end{smallmatrix}\right)\phi_{\alpha_j}\big(\!\left(\begin{smallmatrix}
 T_{00}^j & T_{01}^j \\
 0 & T_{11}^j \\
 \end{smallmatrix}\right)\!\big)\left(\begin{smallmatrix}
 I & X \\
 0 & I \\
 \end{smallmatrix}\right)=\left(\begin{smallmatrix}
 \phi_{\alpha_j}(T_{00}^j) & 0 \\
 0 & \phi_{\alpha_j}(T_{11}^j) \\
 \end{smallmatrix}\right), \phi_{\alpha_j}\in M\ddot{o}b, 1\leq j\leq m.$$
Let $\phi_\alpha=(\phi_{\alpha_1}, \phi_{\alpha_2}, \cdots,\phi_{\alpha_m})$. Then $\phi_\alpha\in M\ddot{o}b^m$ and $\left(\begin{smallmatrix}
 I & -X \\
 0 & I \\
 \end{smallmatrix}\right)\phi_{\alpha}(\widetilde{\textbf{T}})\left(\begin{smallmatrix}
 I & X \\
 0 & I \\
 \end{smallmatrix}\right)=\phi_{\alpha}(\textbf{T}_{00}\oplus\textbf{T}_{11})$.
Since $\textbf{T}_{00},\textbf{T}_{11}$ are both weakly homogeneous, it follows that there exists invertible operator $Y_{\alpha}$ depending on $\alpha$, such that
$Y_{\alpha}^{-1}\left(\begin{smallmatrix}
 I & -X \\
 0 & I \\
 \end{smallmatrix}\right)\phi_{\alpha}(\widetilde{\textbf{T}})\left(\begin{smallmatrix}
 I & X \\
 0 & I \\
 \end{smallmatrix}\right)Y_{\alpha}=\textbf{T}_{00}\oplus\textbf{T}_{11}$.
By using Proposition \ref{3.5pro} again, we obtain that $\left(\begin{smallmatrix}
 I & X \\
 0 & I \\
 \end{smallmatrix}\right)Y_{\alpha}^{-1}\left(\begin{smallmatrix}
 I & -X \\
 0 & I \\
 \end{smallmatrix}\right)\phi_\alpha(\widetilde{\textbf{T}})\left(\begin{smallmatrix}
 I & X \\
 0 & I \\
 \end{smallmatrix}\right)Y_{\alpha}\left(\begin{smallmatrix}
 I & -X \\
 0 & I \\
 \end{smallmatrix}\right)=\widetilde{\textbf{T}}.$
Hence, $\widetilde{\textbf{T}}$ is weakly homogeneous.
\end{proof}

\begin{prop}\label{22217.1}
Let $\text{\bf{T}}_{ii}\in\mathcal{B}^{m}_{n_i}(\mathbb{D}^m),i=0,1$ and $\text{\bf{T}}_{01}\in ran\sigma_{\text{\bf{T}}_{00},\text{\bf{T}}_{11}}$.
Suppose that
$\widetilde{\text{\bf{T}}}=(T_{1}, \cdots, T_{m})\in N\mathcal{FB}_{n_0+n_1}^{m}(\mathbb{D}^m)$ with $T_{i}=\left(\begin{smallmatrix}
 T_{00}^j & T_{01}^j \\
 0 & T_{11}^j \\
 \end{smallmatrix}\right),1\leq j\leq m.$
If $\widetilde{\text{\bf{T}}}$ is weakly homogeneous with respect to $M\ddot{o}b^m$, then $\text{\bf{T}}_{00}\oplus\text{\bf{T}}_{11}$ is also $M\ddot{o}b^m$-weakly homogeneous.
\end{prop}
\begin{proof}
Suppose that there exists $X$ such that $\textbf{T}_{01}=\sigma_{\textbf{T}_{00},\textbf{T}_{11}}(-X)$.
If $\widetilde{\textbf{T}}$ is weakly homogeneous with respect to $M\ddot{o}b^m$, then there exists invertible operator $Y_{\alpha}$ depending on $\alpha=(\alpha_1,\alpha_2,\cdots,\alpha_m)\in \mathbb{D}^m$ and $\phi_\alpha=(\phi_{\alpha_1},\phi_{\alpha_2},\cdots,\phi_{\alpha_m})\in M\ddot{o}b^m$, such that $Y_{\alpha}^{-1}\phi_\alpha(\widetilde{\textbf{T}})Y_{\alpha}=\widetilde{\textbf{T}}. $
By using the strong reducibility of $\widetilde{\textbf{T}}$ in Proposition \ref{3.5pro}, we have $$\left(\begin{smallmatrix}
 I & -X \\
 0 & I \\
 \end{smallmatrix}\right)Y_{\alpha}^{-1}\left(\begin{smallmatrix}
 I & X \\
 0 & I \\
 \end{smallmatrix}\right)\phi_\alpha(\textbf{T}_{00}\oplus\textbf{T}_{11})\left(\begin{smallmatrix}
 I & -X \\
 0 & I \\
 \end{smallmatrix}\right)Y_{\alpha}\left(\begin{smallmatrix}
 I & X \\
 0 & I \\
 \end{smallmatrix}\right)=\textbf{T}_{00}\oplus\textbf{T}_{11}. $$Thus, $\textbf{T}_{00}\oplus\textbf{T}_{11}$ is weakly homogeneous with respect to $M\ddot{o}b^m$.
\end{proof}

\begin{prop}
If $\left(\begin{smallmatrix}
 T_{00} & T_{01} \\
 0 & T_{11} \\
 \end{smallmatrix}\right)\in N\mathcal{FB}_{2}^1(\mathbb{D})$ is weakly homogeneous and $ker \sigma_{\phi_{\alpha}(T_{00}),T_{11}}=\{0\}$ for any $\phi_{\alpha}\in M\ddot{o}b$. Then
$T_{00},T_{11}$ are both weakly homogeneous.
\end{prop}
\begin{proof}
Suppose that there exists $X$ such that $T_{01}=\sigma_{T_{00},T_{11}}(-X)$.
If $\left(\begin{smallmatrix}
 T_{00} & T_{01} \\
 0 & T_{11} \\
 \end{smallmatrix}\right)$ is weakly homogeneous, by Proposition \ref{22217.1}, we see that $T_{00}\oplus T_{11}$ is weakly homogeneous, that is,
$\phi_\alpha(T_{00})\oplus \phi_\alpha(T_{11})$ is similar to $T_{00}\oplus T_{11}$ for any $\phi_{\alpha}\in M\ddot{o}b$.
Note that $\phi_\alpha(T_{ii})\in \mathcal{B}_1^1(\mathbb{D}),i=0,1$,
by the main theorem of \cite{JGJ}, we know that every Cowen-Douglas operator has a unique strongly irreducible decomposition up to similarity, then either $\phi_\alpha(T_{ii})\sim_{s}T_{ii},i=0,1$ or
$\phi_\alpha(T_{00})\sim_{s}T_{11},\phi_\alpha(T_{11})\sim_{s}T_{00}$ for $\phi_\alpha \in M\ddot{o}b$. Since $\ker \sigma_{\phi_{\alpha}(T_{00}),T_{11}}=\{0\}$, we have $\phi_\alpha(T_{ii})\sim_{s}T_{ii},i=0,1$ for any $\phi_\alpha \in M\ddot{o}b$. Hence, $T_{00},T_{11}$ are both weakly homogeneous.
\end{proof}

\begin{prop}\label{2111.27}
Let $T\in\mathcal{B}^1_1(\mathbb{D})$. If $T$ is weakly homogeneous, then there exists $\Psi(\cdot,\cdot):\mathbb{D}\times\mathbb{D}\rightarrow \mathbb{R}^{+}$ such that
$\mathcal{K}_{\Psi}(\alpha,\phi_\alpha(w))|\phi'_\alpha(w)|^2+\mathcal{K}_{\Psi}(\alpha,w)=0,$
where $\mathcal{K}_{\Psi}(\alpha,w)=-\frac{\partial^{2}}{\partial w\partial \bar{w}}\log\Psi(\alpha,w)$ and $\phi_\alpha \in M\ddot{o}b$. In particular, $\Psi(\alpha,w)$ also satisfies
$\mathcal{K}_{\Psi}(w,w)=-\frac{\mathcal{K}_{\Psi}(w,0)}{(1-|w|^2)^2}$ and $\mathcal{K}_{\Psi}(0,-w)=-\mathcal{K}_{\Psi}(0,w)$.
\end{prop}
\begin{proof}
Without losing generality, we assume that $\phi_\alpha(w)=\frac{\alpha-w}{1-\bar{\alpha}w},\alpha,w\in \mathbb{D}$. If $T$ is weakly homogeneous, then for any $\alpha\in \mathbb{D}$, $T$ is similar to $\phi_\alpha(T)$. Let $e$ be a non-zero section of $E_T$ associated with $T$.
Note that $\phi_\alpha(T)\in\mathcal{B}^1_1(\mathbb{D})$ and $e(\phi_\alpha(w))\in\ker(\phi_\alpha(T)-w)$.
By Theorem \ref{main9.26}, we find that a bounded operator $X_\alpha$ and $\psi_\alpha\in H^\infty(\mathbb{D})$ depending on $\alpha$,
such that $\|e(\phi_\alpha(w))\|^2=|\psi_\alpha(w)|^2(\|e(w)\|^2+\|X_\alpha(e(w))\|^2),w\in \mathbb{D}.$
Further, we have $\frac{\partial^{2}}{\partial w\partial \bar{w}}\log\|e(\phi_\alpha(w))\|^2=\frac{\partial^{2}}{\partial w\partial \bar{w}}\log\|e(w)\|^2+\frac{\partial^{2}}{\partial w\partial \bar{w}}\log\Big(1+\frac{\|X_\alpha(e(w))\|^2}{\|e(w)\|^2}\Big),w\in \mathbb{D}.$
Define $\Psi(\cdot,\cdot):\mathbb{D}\times\mathbb{D}\rightarrow \mathbb{R}^{+}$ as $\Psi(\alpha,w):=1+\frac{\|X_\alpha(e(w))\|^2}{\|e(w)\|^2}$ and $\mathcal{K}_{\Psi}(\alpha,w)=-\frac{\partial^{2}}{\partial w\partial \bar{w}}\log\Psi(\alpha,w)$. Thus we infer
\begin{equation}\label{211141}
\mathcal{K}_{\phi_\alpha(T)}(w)=\mathcal{K}_{T}(w)+\mathcal{K}_{\Psi}(\alpha,w).
\end{equation}
From the arbitrariness of $w\in \mathbb{D}$ in equation (\ref{211141}), we replace $w$ with $\phi_\alpha(w)$. It follows that
\begin{equation}\label{211142}
\mathcal{K}_{\phi_\alpha(T)}(\phi_\alpha(w))=\mathcal{K}_{T}(\phi_\alpha(w))+\mathcal{K}_{\Psi}(\alpha,\phi_\alpha(w)).
\end{equation}
It can be obtained by a simple application of the chain rule that
$\mathcal{K}_{T}(w)=\mathcal{K}_{\phi_\alpha(T)}(\phi_\alpha(w))|\phi'_\alpha(w)|^{2}. $
Then equation (\ref{211142}) can be transformed into
\begin{equation}\label{211143}
\mathcal{K}_{T}(w)|\phi'_\alpha(w)|^{-2}=\mathcal{K}_{T}(\phi_\alpha(w))+\mathcal{K}_{\Psi}(\alpha,\phi_\alpha(w)).
\end{equation}
Using the chain rule again, we have
$\mathcal{K}_{\phi_\alpha(T)}(w)=-\frac{\partial^{2}}{\partial w\partial \bar{w}}\log\|e(\phi_\alpha(w))\|^2=\mathcal{K}_{T}(\phi_\alpha(w))|\phi'_\alpha(w)|^{2}.$
Then equation (\ref{211141}) is equivalent to
\begin{equation}\label{211144}
\mathcal{K}_{T}(\phi_\alpha(w))=\mathcal{K}_{T}(w)|\phi'_\alpha(w)|^{-2}+\mathcal{K}_{\Psi}(\alpha,w)|\phi'_\alpha(w)|^{-2}.
\end{equation}
Combining equations (\ref{211143}) and (\ref{211144}), we obtain that
\begin{equation}\label{211145}
\mathcal{K}_{\Psi}(\alpha,\phi_\alpha(w))|\phi'_\alpha(w)|^{2}+\mathcal{K}_{\Psi}(\alpha,w)=0.
\end{equation}

Note that $\phi_{\alpha}'(w)=\frac{|\alpha|^{2}-1}{(1-\bar{\alpha}w)^2}$ and $\phi_\alpha(0)=\alpha$. It follows that
$\mathcal{K}_{\Psi}(\alpha,\alpha)(1-|\alpha|^{2})^{2}+\mathcal{K}_{\Psi}(\alpha,0)=0$ for all $\alpha\in\mathbb{D}$, that is,
$\mathcal{K}_{\Psi}(w,w)=-\frac{\mathcal{K}_{\Psi}(w,0)}{(1-|w|^2)^2}$ for all $w\in\mathbb{D}$.
Since $\phi_0(w)=-w$, we imply $\mathcal{K}_{\Psi}(0,-w)+\mathcal{K}_{\Psi}(0,w)=0$.
\end{proof}

\begin{prop}
Let $T\in\mathcal{B}^1_1(\mathbb{D})$ and $e$ be a non-zero section of $E$ determined by $T$.
For any $\phi_\alpha\in M\ddot{o}b$, $E_\alpha$ is the Hermitian holomorphic vector bundle associated with $\phi_\alpha(T)$.
If $T$ is weakly homogeneous, then there exists vector bundle $F^{\alpha}$ with $F^{\alpha}(w)=\bigvee\{\left(\begin{smallmatrix}
 I \\
 X_\alpha \\
 \end{smallmatrix}\right)e(w)\}$, such that $E\otimes E_{\alpha}\sim_u F^\alpha\otimes F_{\alpha}^\alpha$, where $X_\alpha$ is a bounded operator depending on $\alpha$ and $F_{\alpha}^\alpha(w)=\bigvee\{\left(\begin{smallmatrix}
 I \\
 X_\alpha \\
 \end{smallmatrix}\right)e(\phi_\alpha(w))\},w\in\mathbb{D}$.
\end{prop}
\begin{proof}
Without losing generality, we assume that $\phi_\alpha(w)=\frac{\alpha-w}{1-\bar{\alpha}w}$. If $T$ is weakly homogeneous, then $T\sim_s\phi_\alpha(T), \alpha\in \mathbb{D}$.
It is easy to see that $\phi_\alpha(T)\in\mathcal{B}^1_1(\mathbb{D})$ and $e(\phi_\alpha(w))\in\ker(\phi_\alpha(T)-w)$.
By Theorem \ref{main9.26}, we know that there exists a bounded operator $X_\alpha$ and $\psi_\alpha\in H^\infty(\mathbb{D})$ depending on $\alpha$,
such that $\|e(\phi_\alpha(w))\|^2=|\psi_\alpha(w)|^2(\|e(w)\|^2+\|X_\alpha(e(w))\|^2),w\in \mathbb{D}.$
This is equivalent to
$\frac{\|e(\phi_\alpha(w))\|^2}{\|e(w)\|^2}=|\psi_\alpha(w)|^2\big(1+\frac{\|X_\alpha(e(w))\|^2}{\|e(w)\|^2}\big)$ and
$\frac{\|e(w)\|^2}{\|e(\phi_\alpha(w))\|^2}=|\psi_\alpha(\phi_\alpha(w))|^2\big(1+\frac{\|X_\alpha(e(\phi_\alpha(w)))\|^2}{\|e(\phi_\alpha(w))\|^2}\big)$, since $\phi_\alpha(\phi_\alpha(w))=1,\alpha,w\in \mathbb{D}$.
Then
$$\begin{array}{lll}
&&|\psi_\alpha(w)\psi_\alpha(\phi_\alpha(w))|^2\cdot \frac{(\|e(w)\|^2+\|X_\alpha(e(w))\|^2)(\|e(\phi_\alpha(w))\|^2+\|X_\alpha(e(\phi_\alpha(w)))\|^2)}{\|e(w)\|^2\|e(\phi_\alpha(w))\|^2}\\
&=&|\psi_\alpha(w)\psi_\alpha(\phi_\alpha(w))|^2\cdot \frac{\Big\|\left(\begin{smallmatrix}
 e(w) \\
 X_\alpha(e(w)) \\
 \end{smallmatrix}\right)\otimes\left(\begin{smallmatrix}
 e(\phi_\alpha(w)) \\
 X_\alpha(e(\phi_\alpha(w))) \\
 \end{smallmatrix}\right)\Big\|^2}{\|e(w)\otimes e(\phi_\alpha(w))\|^2}\\
&=&1, \alpha,w\in \mathbb{D}.
\end{array}$$
Further, we have
\begin{equation}\label{2202.23}
\frac{\partial^{2}}{\partial w\partial \bar{w}}\log\frac{\big\|\left(\begin{smallmatrix}
 e(w) \\
 X_\alpha(e(w)) \\
 \end{smallmatrix}\right)\otimes\left(\begin{smallmatrix}
 e(\phi_\alpha(w)) \\
 X_\alpha(e(\phi_\alpha(w))) \\
 \end{smallmatrix}\right)\big\|^2}{\|e(w)\otimes e(\phi_\alpha(w))\|^2}=0,\alpha,w\in \mathbb{D}.
\end{equation}
Let $F^\alpha(w)=\bigvee\{\left(\begin{smallmatrix}
 I \\
 X_\alpha \\
 \end{smallmatrix}\right)e(w)\}$ and $F_{\alpha}^\alpha(w)=\bigvee\{\left(\begin{smallmatrix}
 I \\
 X_\alpha \\
 \end{smallmatrix}\right)e(\phi_\alpha(w))\}$ for $\alpha,w\in \mathbb{D}$. By Theorem \ref{3.5you} and equation (\ref{2202.23}), we infer $E\otimes E_{\alpha}\sim_u F^\alpha\otimes F_{\alpha}^\alpha$.
\end{proof}

\begin{rem}
Let $T$ be a Cowen-Douglas operator with index one. By Proposition 3.2 in \cite{HJZ}, the first and second authors joint with L. Zhao proved that for any $\phi_\alpha \in M\ddot{o}b$,
if the holomorphic Hermitian
vector bundle $E_\alpha$ associated with $\phi_\alpha(T)$ is congruent to $E_T\otimes\mathcal{L}_\alpha$ for some line bundle $\mathcal{L}_\alpha$,
then $T$ is homogeneous.
In this case, when the index of $T$ is one, the $\mathcal{K}_{\Psi}$ in Proposition \ref{2111.27} is zero, since $\|e(\phi_\alpha(w))\|^2=|\varphi_\alpha(w)|^2\|e(w)\|^2$ for some holomorphic function $\varphi$ and $\alpha,w\in\mathbb{D}$.
But
if $\|e(\phi_\alpha(w))\|^2=(1+|\varphi_\alpha(w)|^2)\|e(w)\|^2$, then $T$ is not a homogeneous operator, since $1+|\varphi_\alpha(w)|^2$ is not the square of the Modulus of some holomorphic function.
Although it can be regarded as $X_\alpha e(w)=\varphi_\alpha(w)e(w)$, the equation $\mathcal{K}_{\Psi}(\alpha,\phi_\alpha(w))|\phi'_\alpha(w)|^2+\mathcal{K}_{\Psi}(\alpha,w)=0$,
a necessary condition of $T$ to be a weakly homogeneous operator, does not hold. We have
$\mathcal{K}_{\Psi}(\alpha,w)=-\frac{|\varphi_\alpha'(w)|^2}{(1+|\varphi_\alpha(w)|^2)^2}$ and
$\mathcal{K}_{\Psi}(\alpha,\phi_\alpha(w))=-\frac{|\varphi_\alpha'(\phi_\alpha(w))|^2}{(1+|\varphi_\alpha(\phi_\alpha(w))|^2)^2}$.
The reason for this phenomenon may be
$(I+X^*_\alpha X_\alpha)^{\frac{1}{2}}$ is not the intertwining of $T$ and $\phi_\alpha(T)$.
\end{rem}

In what follows, we assume that Hilbert space $\mathcal{H}_{i},i\geq0$ is analytical function space with reproducing kernel $K^i(z,w)$, where
$K^0(z,w)=-\frac{\ln(1-z\bar{w})}{z\bar{w}}$, $K^n(z,w)=\frac{1}{(1-z\bar{w})^n}, n\geq1, z,w\in\mathbb{D}$.
Let $T$ be a Cowen-Douglas operator and $T\sim_{u}(M_z^*,\mathcal{H}_K)$. When $T \in \mathcal{B}_n^1(\mathbb{D})$ and is contractive,
M. Uchiyama in \cite{MU} provide a necessary and sufficient condition for $T$ is similar to the $n$ times copies of $M^*_z$ on Hardy space, which is that
there exist positive constants $m,M$ such that $m\sum_{i=1}^n|x_i|^2\leq(1-|w|^2)\langle K(\bar{w}, \bar{w})\xi,\xi\rangle\leq M\sum_{i=1}^n|x_i|^2$
for any $w\in\mathbb{D}$ and $\xi=\sum_{i=1}^nx_i\xi_i,x_i\in\mathbb{C},\xi_i=(0,\cdots,0,1,0,\cdots0)^T$ with 1 on the $i$th position.
When $T \in \mathcal{B}_1^1(\mathbb{D})$ and is $n$-hypercontractive, the second and third authors in \cite{jiandji} show that $T$ is similar to $M_z^*$ on $(\mathcal{H}_n,K^n)$ if and only if $\frac{K(w,w)}{K^n(w,w)}$ is bounded and bounded below from zero.
For each $n\geq1$, we know that the multiplication operator on $(\mathcal{H}_n,K^n)$ is homogeneous in \cite{Misra} given by G. Misra.
It is well known that an operator is similar to a homogeneous is a weakly homogeneous. The $n-$hypercontraction $T$, $n\geq1$ determined by the similarity above is weakly homogeneous.
For some positive definite kernels $K$, it is shown that in Theorem 5.3 of \cite{Ghara} the
multiplication operator $M_z$ on $(\mathcal{H},KK^{n}), n>0$, is a weakly
homogeneous operator due to S. Ghara.

In the following, we provide some methods to obtain weakly homogeneous tuples.
We let $(\mathbf{M}_z,\mathcal{H}_{K})$ denote the tuple of multiplication operators acting on Hilbert space $\mathcal{H}_{K}$ determined by the unique non-negative definite kernel $K$.

\begin{prop}\label{220419.1}
Let $\mathcal{H}_{K_i}$ be the analytic function space with reproducing kernel $K_i$ over $\Omega\subset\mathbb{C}^m$, $i=0,1$.
Suppose that the tuple of multiplication operators $\mathbf{M}_z=(M_{z_1},\cdots,M_{z_m})$ acting on $\mathcal{H}_{K_0}$ and $\mathcal{H}_{K_1}$ are bounded.
If the identity mapping $id:\mathcal{H}_{K_0}\rightarrow\mathcal{H}_{K_1}$ is bounded, then
$(\mathbf{M}_z,\mathcal{H}_{K_0+K_1})$ is similar to $(\mathbf{M}_z,\mathcal{H}_{K_1})$.
\end{prop}
\begin{proof}
For any analytic function $f$ and $\xi\in\mathbb{C}^n$, we have $$\langle f,(id)^*K_1(\cdot,w)\xi\rangle_{\mathcal{H}_{K_0}}=\langle f,K_1(\cdot,w)\xi\rangle_{\mathcal{H}_{K_1}}=\langle f(w),\xi\rangle_{\mathbb{C}^n}=\langle f,K_0(\cdot,w)\xi\rangle_{\mathcal{H}_{K_0}},w\in\Omega.$$
By the arbitrariness of $f$, it follows that $(id)^*K_1(\cdot,w)\xi=K_0(\cdot,w)\xi$ for all $w\in\Omega$ and $\xi\in\mathbb{C}^n$.
Let $X:=(id)^*:\mathcal{H}_{K_1}\rightarrow\mathcal{H}_{K_0}$. If $id:\mathcal{H}_{K_0}\rightarrow\mathcal{H}_{K_1}$ is bounded, then so is $X$.
Since $\mathcal{H}_{K_i}$ is the analytic function space with reproducing kernel $K_i$, $i=0,1$, then so is $\mathcal{H}_{K_0+K_1}$ with $K_0+K_1$ in \cite{CS2,Salinas} due to R. E. Curto and N. Salinas.
We know that the tuple of multiplication operators $\mathbf{M}_z=(M_{z_1},\cdots,M_{z_m})$ satisfies $$\bigvee_{w\in \Omega^*}\mbox{ker}(\mathbf{M}_z^*-w)=\bigvee_{w\in \Omega^*}\{K_i(\cdot,\bar{w})\xi,\xi\in\mathbb{C}^n\}=\mathcal{H}_{K_i},\,i=0,1,$$ then the reproducing kernel of $\mathbf{M}_z$ on $\mathcal{H}_{K_0+K_1}$ is
$$\begin{array}{lll}
&&\big(\langle (K_0(\cdot,w)+K_1(\cdot,w))e_j,(K_0(\cdot,w)+K_1(\cdot,w))e_i\rangle_{\mathcal{H}_{K_0+K_1}}\big)_{i,j=1}^n\\
&=&\big(\langle (K_0(w,w)+K_1(w,w))e_j,e_i\rangle_{\mathcal{H}_{K_0+K_1}}\big)_{i,j=1}^n\\
&=&\big(\langle K_0(w,w)e_j,e_i\rangle_{\mathcal{H}_{K_0}}\big)_{i,j=1}^n+\big(\langle K_1(w,w)e_j,e_i\rangle_{\mathcal{H}_{K_1}}\big)_{i,j=1}^n
\end{array}$$
for any orthonormal basis $\{e_i\}_{i=1}^n$ of $\mathbb{C}^n$. Note that
$$\begin{array}{lll}
&&\big(\langle (I_{\mathcal{H}_{K_1}}+X^*X)^{\frac{1}{2}}K_1(\cdot,w)e_j,(I_{\mathcal{H}_{K_1}}+X^*X)^{\frac{1}{2}}K_1(\cdot,w)e_i\rangle_{\mathcal{H}_{K_1}}\big)_{i,j=1}^n\\
&=&\big(\langle (I_{\mathcal{H}_{K_1}}+X^*X)K_1(\cdot,w)e_j,K_1(\cdot,w)e_i\rangle_{\mathcal{H}_{K_1}}\big)_{i,j=1}^n\\
&=&\big(\langle K_1(\cdot,w)e_j,K_1(\cdot,w)e_i\rangle_{\mathcal{H}_{K_1}}\big)_{i,j=1}^n+\big(\langle XK_1(\cdot,w)e_j,XK_1(\cdot,w)e_i\rangle_{\mathcal{H}_{K_0}}\big)_{i,j=1}^n\\
&=&\big(\langle K_1(\cdot,w)e_j,K_1(\cdot,w)e_i\rangle_{\mathcal{H}_{K_1}}\big)_{i,j=1}^n+\big(\langle K_0(\cdot,w)e_j,K_0(\cdot,w)e_i\rangle_{\mathcal{H}_{K_0}}\big)_{i,j=1}^n,
\end{array}$$
then means that the reproducing kernels of $(I_{\mathcal{H}_{K_1}}+X^*X)^{\frac{1}{2}}\mathbf{M}_z(I_{\mathcal{H}_{K_1}}+X^*X)^{-\frac{1}{2}}$ on $\mathcal{H}_{K_1}$ and $\mathbf{M}_z$ on $\mathcal{H}_{K_0+K_1}$ are the same.
By Remark 3.8 of \cite{CS2}, we have $(I_{\mathcal{H}_{K_1}}+X^*X)^{\frac{1}{2}}\mathbf{M}_z$ on $\mathcal{H}_{K_1}$ is unitarily equivalent to $\mathbf{M}_z$ on $\mathcal{H}_{K_0+K_1}$.
Since $(I_{\mathcal{H}_{K_1}}+X^*X)^{\frac{1}{2}}$ is invertible, then $(\mathbf{M}_z,\mathcal{H}_{K_0+K_1})$ is similar to $(\mathbf{M}_z,\mathcal{H}_{K_1})$.
\end{proof}
Let $K_0,K_1:\Omega\times\Omega\rightarrow \mathcal{M}_n(\mathbb{C})$. We write $K_0\succeq0$ means that $K_0$ is a non-negative definite kernel.
We write $K_0\preceq K_1$ or $K_1\succeq K_0$, if $K_0$ and $K_1$ are two non-negative kernels satisfying $K_1-K_0\succeq0$ on $\Omega\times\Omega$.
\begin{prop}\label{220507.1}
Let $\mathcal{H}_{K_i}$ be Hilbert space determined by reproducing kernel $K_i$ over $\Omega\subset\mathbb{C}^m$, $i=0,1$.
Suppose that the multiplication operator $\mathbf{M}_{z}=(M_{z_1},\cdots,M_{z_m})$ is bounded on $(\mathcal{H}_{K_i} ,K_i)$ for $i=0,1$ and $K_1\succeq K_0$ on $\Omega\times\Omega$.
Then $(\mathbf{M}_z,\mathcal{H}_{K_0+K_1})$ is similar to $(\mathbf{M}_z,\mathcal{H}_{K_1})$.
\end{prop}
\begin{proof}
Recall that if $K_1\succeq K_0$, then
$(\mathcal{H}_{K_0},K_0)\subset(\mathcal{H}_{K_1},K_1)$ and $\|h\|_{\mathcal{H}_{K_1}}\leq\|h\|_{\mathcal{H}_{K_0}}$ for $h\in(\mathcal{H}_{K_0},K_0)$ in Theorem 6.25 of \cite{PR58}.
Let $id:\mathcal{H}_{K_0}\rightarrow\mathcal{H}_{K_1}$ be the identity mapping. For any $h\in(\mathcal{H}_{K_0},K_0)$, we have $\frac{\|id(h)\|_{\mathcal{H}_{K_1}}}{\|h\|_{\mathcal{H}_{K_0}}}=\frac{\|h\|_{\mathcal{H}_{K_1}}}{\|h\|_{\mathcal{H}_{K_0}}}\leq1$. It follows that
$id:\mathcal{H}_{K_0}\rightarrow\mathcal{H}_{K_1}$ is bounded.
By Proposition \ref{220419.1}, we obtain that $(\mathbf{M}_z,\mathcal{H}_{K_0+K_1})$ is similar to $(\mathbf{M}_z,\mathcal{H}_{K_1})$.
\end{proof}
In \cite{AM}, A. Kor\'{a}nyi and G. Misra explicitly construct all homogeneous holomorphic Hermitian vector bundles on $\mathbb{D}$ and give the forms of all operators similar to homogeneous operators.
\begin{lem}\cite{AM}\label{5505.1301}
Any homogeneous operator in the Cowen-Douglas class is similar to the direct sum (with multiplicity) of weighted block shifts, which is the direct sum of an ordinary
(unweighted) block shift and a Hilbert-Schmidt operator.
\end{lem}

\begin{rem}\label{220508.3}
If we choose $K_1$ as an arbitrary homogeneous kernel (the tuple $\mathbf{M}_{z}$ of multiplications acting on $(\mathcal{H}_{K_1},K_1)$ is homogeneous), $K_0$ is an arbitrary nonnegative definite kernel and satisfies $K_1\succeq K_0$, then $(\mathbf{M}_z,\mathcal{H}_{K_0+K_1})$ is weakly homogeneous by Proposition \ref{220507.1}.
In the following, we simply lists some examples, but does not give detailed proof.

For any nonnegative integer $k$ and $\{i_0,\cdots,i_l\}\subset\{1,\cdots,k\}$,
$(M_z,\mathcal{H}_{\sum_{j=0}^lK^{i_j}})$ is similar to $(M_z,\mathcal{H}_{K^r})$, where $r=\max\{i_0,\cdots,i_l\}$.
Thus, $(M_z,\mathcal{H}_{\sum_{j=0}^lK^{i_j}})$ is weakly homogeneous. From Lemma \ref{5505.1301}, it is also similar to a weighted block shift.
Let $K^{(\lambda)}(z,w)=\prod_{i=1}^m\frac{1}{(1-z_i\bar{w}_i)^\lambda}$ and $\tilde{K}^{(\alpha)}(z,w)=\prod_{i=1}^m\frac{1}{(1-z_i\bar{w}_i)^{\alpha_i}}$ for
$z=(z_1,\cdots,z_m),w=(w_1,\cdots,w_m)\in\mathbb{D}^m$,
positive integer $\lambda$ and  $\alpha=(\alpha_1,\cdots,\alpha_m)\in \mathbf{Z}_{+}^{m}$.
Then the tuples of multiplications acting on $(\mathcal{H}_{K^{(\lambda)}+K^{(\mu)}},K^{(\lambda)}+K^{(\mu)})$,
$(\mathcal{H}_{\tilde{K}^{(\alpha)}+\tilde{K}^{(\beta)}},\tilde{K}^{(\alpha)}+\tilde{K}^{(\beta)})$ are weakly homogeneous for $0<\lambda\leq\mu$ and $\alpha\leq\beta\in \mathbf{Z}_{+}^{m}$.

Next, we consider an example of homogeneous tuple on $\mathbb{B}^m$, the open unit ball $\{w:|w|<1\}$ in $\mathbb{C}^m$.
Let $\mathcal{H}_{K_{(n)}}$ be the reproducing kernel
Hilbert space determined by the kernel function
$K_{(n)}(z,w)=\frac{1}{(1-\langle z,w\rangle)^n}=\sum_{\alpha\in \mathbf{Z}_{+}^{m}}\frac{(n+|\alpha|-1)!}{\alpha!(n-1)!}z^\alpha\bar{w}^\alpha$ for positive integer $n$,
$z=(z_1,\cdots,z_m),w=(w_1,\cdots,w_m)\in\mathbb{B}^m$ and $\langle z,w\rangle=z_1\bar{w}_1+\cdots+z_m\bar{w}_m$.
The tuples of multiplications on $\mathcal{H}_{K_{(n)}}$ is homogeneous respect to $G$, the bi-holomorphic
automorphism group of the ball $\mathbb{B}^m$.
Then the tuples of multiplications acting on $(\mathcal{H}_{K_{(n_1)}+K_{(n_2)}},K_{(n_1)}+K_{(n_2)})$ is weakly homogeneous for $n_1,n_2>0$.

Next, we consider the operator class $\mathcal{FB}_n^1(\Omega)$, which has a flag structure and introduced in \cite{JJKMCR,JJKM}. Let $T\in\mathcal{FB}_n^1(\Omega)$ with a form in Definition 3.1 of \cite{JJKM} and $T_i\sim_{u}(M_z^*,\mathcal{H}_{\tilde{K}_i}),0\leq i\leq n-1$. If $S_{i,i+1}\tilde{K}_{i+1}(\cdot,\bar{w})=\tilde{K}_{i}(\cdot,\bar{w})$ and it is not zero on $w\in\Omega_0\subset\Omega$, then
$(M_z,\mathcal{H}_{\sum_{j=0}^lK_{j}})$ is similar to $(M_z,\mathcal{H}_{K_l})$ on $\Omega_0$.
Naturally, this result can be extended to some tuples in the class $\mathcal{FB}_2^m(\Omega)$.

Let $\mathcal{H}^{(\lambda,\mu)}$ be the reproducing kernel Hilbert space with the kernel $K^{(\lambda,\mu)}(z,w)=\left(\begin{smallmatrix}
\frac{1}{(1-z\bar{w})^\lambda} & \frac{z}{(1-z\bar{w})^{\lambda+1}} \\
\frac{\bar{w}}{(1-z\bar{w})^{\lambda+1}} & \frac{\frac{1}{\lambda}+\mu+z\bar{w}}{(1-z\bar{w})^{\lambda+2}}\\
\end{smallmatrix}\right)$ on $\mathbb{D}\times\mathbb{D}$, $\lambda,\mu>0$.
It is shown in \cite{KM1,Wilkins} that every irreducible homogeneous operators in $\mathcal{B}_2^1(\mathbb{D})$ must be unitarily equivalent to
the adjoint of multiplication operator acting on $\mathcal{H}^{(\lambda,\mu)}$
for some $\lambda,\mu>0$. Let $\lambda_i,\mu_i>0,i=0,1$. By a direct computation, it can be seen $K^{(\lambda_1,\mu_1)}\succeq K^{(\lambda_0,\mu_0)}$ when $\lambda_1\geq\lambda_0$ and $\mu_1\geq\frac{1}{\lambda_0}-\frac{1}{\lambda_1}+\mu_0$.
Thus, from Proposition \ref{220507.1}, we obtain that $(M_z,\mathcal{H}_{K^{(\lambda_0,\mu_0)}+K^{(\lambda_1,\mu_1)}})$ is similar to $(M_z,\mathcal{H}_{K^{(\lambda_1,\mu_1)}})$ and $(M_z,\mathcal{H}_{K^{(\lambda_0,\mu_0)}+K^{(\lambda_1,\mu_1)}})$ is weakly homogeneous.
\end{rem}

Let $\tilde{\mathcal{H}}^{(\lambda,\mu)}=\mathcal{H}_{\lambda_1}\otimes\cdots\otimes\mathcal{H}_{\lambda_{m-1}}\otimes\mathcal{H}^{(\lambda_m,\mu)}$, $\lambda_i,\mu>0, 1\leq i\leq m$.
Then $\tilde{\mathcal{H}}^{(\lambda,\mu)}$ is the reproducing kernel Hilbert space determined by the kernel $$\tilde{K}^{(\lambda, \mu)}(z,w)=\Big(\prod_{i=1}^{m-1}K^{\lambda_i}(z_i,w_i)\Big)K^{(\lambda_{m},\mu)}(z_m,w_m)$$ on $\mathbb{D}^m\times\mathbb{D}^m$ for a tuple of positive real numbers $\lambda=(\lambda_1, \cdots,\lambda_m)$ and $\mu>0$.
Based on the above discussion of kernels of $\tilde{K}^{(\alpha)}$ and $K^{(\lambda,\mu)}$, for $\lambda,\tilde{\lambda}\in \mathbf{Z}_{+}^{m}$ and $\mu,\tilde{\mu}>0$, if
$\lambda\geq\tilde{\lambda}$ and $\mu\geq\frac{1}{\tilde{\lambda}_m}-\frac{1}{\lambda_m}+\tilde{\mu}$, we obtain that
$\tilde{K}^{(\lambda, \mu)}\succeq\tilde{K}^{(\tilde{\lambda}, \tilde{\mu})}$ on $\mathbb{D}^m\times\mathbb{D}^m$.
It follows from Proposition \ref{220507.1} that $(\mathbf{M}_z,\mathcal{H}_{\tilde{K}^{(\lambda, \mu)}+\tilde{K}^{(\tilde{\lambda}, \tilde{\mu})}})$ is similar to $(\mathbf{M}_z,\mathcal{H}_{\tilde{K}^{(\lambda, \mu)}})$.
By Theorem 6.7 of \cite{DH}, we know that each homogeneous tuple of operators in $\mathcal{B}_2^m(\mathbb{D}^m)$ with respect to M$\ddot{o}$b$^m$ is unitarily equivalent to the adjoint
of the tuple of multiplication operators on the reproducing kernel Hilbert space $\tilde{\mathcal{H}}^{(\lambda,\mu)}$ with the kernel
$\tilde{K}^{(\lambda, \mu)}$, where $\lambda=(\lambda_1, \cdots,\lambda_m)$ is a tuple of positive real numbers and $\mu>0$.
Thus, $(\mathbf{M}_z,\mathcal{H}_{\tilde{K}^{(\lambda, \mu)}+\tilde{K}^{(\tilde{\lambda}, \tilde{\mu})}})$ is weakly homogeneous.

A question raised by K. Zhu in \cite{ZKH2} is that for the multiplication operator $M_{z}$ on the Bergman space and its two invariant subspaces $I$ and $J$, when are the two restriction operators $M_{z}|_I$ and $M_{z}|_J$ are similar?
For this problem, we replace the multiplication operator on Bergman space with the adjoint of the direct sum of multiplication operators on two analytic functions Hilbert spaces and give a sufficient condition.
\begin{rem}
If $\mathcal{H}_{K_0}$ and $\mathcal{H}_{K_1}$ are analytic functions Hilbert spaces with reproducing kernels $K_0$ and $K_1$, respectively, then so is $\mathcal{H}_{K_0+K_1}$ with reproducing kernel $K_0+K_1$ due to \cite{CS2,Salinas}.
In order not to cause misunderstanding, we note that the multiplication operator on $\mathcal{H}_{K_i}$ is $M_{z,i},i=0,1$.
Let $\mathcal{M}=\{(g,-g)\in\mathcal{H}_{K_0}\oplus\mathcal{H}_{K_1}:g\in\mathcal{H}_{K_0}\cap\mathcal{H}_{K_1}\}\subset\mathcal{H}_{K_0}\oplus\mathcal{H}_{K_1}$.
An implicit conclusion in \cite{Salinas} is $(M_z,\mathcal{H}_{K_0+K_1})$ is unitarily equivalent to $P_{\mathcal{M}^\bot}(M_{z,0}\oplus M_{z,1})|_{\mathcal{M}^\bot}$.
For any $(g,-g)\in\mathcal{M}$, we have $(M_{z,0}\oplus M_{z,1})(g,-g)=z(g,-g)\in\mathcal{M}$ and then $\mathcal{M}$ is an invariant subspace of $M_{z,0}\oplus M_{z,1}$.
A simple calculation leads to $\mathcal{M}^\bot$ is an invariant subspace of $M_{z,0}^*\oplus M_{z,1}^*$.
It follows that $(M_z^*,\mathcal{H}_{K_0+K_1})$ is unitarily equivalent to $(M_{z,0}^*\oplus M_{z,1}^*)|_{\mathcal{M}^\bot}$.
In particularly, if $K_1\succeq K_0$, by Proposition \ref{220507.1}, we know that $(M_z^*,\mathcal{H}_{K_0+K_1})$ is similar to $(M_{z,1}^*,\mathcal{H}_{K_1})$.
Note that $(M_{z,1}^*,\mathcal{H}_{K_1})=(M_{z,0}^*\oplus M_{z,1}^*)|_{\mathcal{H}_{K_1}}$.
Hence, $(M_{z,0}^*\oplus M_{z,1}^*)|_{\mathcal{H}_{K_1}}$ is similar to $(M_{z,0}^*\oplus M_{z,1}^*)|_{\mathcal{H}_{\mathcal{M}^\bot}}$.
\end{rem}

\end{document}